\documentclass[11pt]{article}
\newcommand{\version}{October 2007}

\usepackage{amsthm,amsfonts,amsmath,amssymb}
\usepackage{graphicx} 

\swapnumbers

\theoremstyle{plain}
\newtheorem{thm}{THEOREM}[section]
\newtheorem{lm}[thm]{LEMMA}

\newtheorem{prop}[thm]{PROPOSITION}

\newtheorem{defi}[thm]{DEFINITION}

\theoremstyle{remark}

\newcommand{\upchi}{\raise1pt\hbox{$\chi$}}
\newcommand{\R}{{\mathord{\mathbb R}}}

\newcommand{\dd}{\ensuremath \, {\rm d}}

\newcommand{\id}{\mbox{\rm Id}}

\numberwithin{equation}{section} 
\begin{document}

\title{\bf{ Subadditivity of the Entropy and its Relation to Brascamp--Lieb Type Inequalities}}
\author{\vspace{5pt} Eric A. Carlen$^{1}$ and Dario Cordero--Erausquin$^{2}$ \\
\vspace{5pt}\small{$1.$ Department of Mathematics, Rutgers University,}\\[-6pt]
\small{Hill Center,
110 Frelinghuysen Road
Piscataway NJ 08854-8019 USA}\\
\vspace{5pt}\small{$2.$ 
Institut de Math\'ematiques de Jussieu,
Universit\'e Pierre et Marie Curie - Paris 6,} \\[-6pt]
\small{ 4 Place Jussieu, 75252 Paris Cedex 05, France }\\ }
\date{\version}
\maketitle
\footnotetext                                                                         
[1]{Work partially
supported by U.S. National Science Foundation
grant DMS 06-00037.   \\
\copyright\, 2007 by the authors. This paper may be reproduced, in its
entirety, for non-commercial purposes.}

\begin{abstract}

We prove a general duality result showing that a Brascamp--Lieb type inequality is equivalent to 
an inequality expressing subadditivity of the entropy, with a complete correspondence of best
constants and cases of equality. This opens a new approach to the proof of Brascamp--Lieb type inequalities,
via subadditivity of the entropy. We illustrate the utility of this approach by 
proving  a  general inequality expressing the 
subadditivity property of the entropy on $\R^n$, and 
fully determining the cases of equality.  As a consequence of the duality mentioned above, we obtain a
simple new proof of the classical Brascamp--Lieb inequality, and also a fully explicit determination of all of the
cases of equality. We also deduce several other consequences of the general subadditivity inequality, including
a generalization of Hadamard's inequality for determinants.  Finally, we also 
prove a second duality theorem
relating superadditivity  of the Fisher information and 
a sharp convolution type inequality for the fundamental eigenvalues of Schr\"odinger operators.  
Though we focus mainly on the case of random variables in $\R^n$ in this paper, we  discuss
extensions to other settings as well.

\end{abstract}

\bigskip
\centerline{Mathematics Subject Classification Numbers: 26D15, 94A17 }

\newpage

\section{Introduction} \label{intro}

Let $(\Omega,{\cal S},\mu)$ be a measure space, and let $f$ be a probability density on 
$(\Omega,{\cal S},\mu)$. That is, $f$ is a non negative integrable function on $\Omega$
with $\int_\Omega f{\rm d}\mu = 1$.   On the convex subset of probability densities 
\begin{equation}\label{finen}
\left\{\ f\ :\ \int_\Omega  f\ln (1+f)\dd \mu < \infty\ \right\}\ ,
\end{equation}
the {\em entropy} of $f$, $S(f)$, is defined by
$$S(f) = \int_\Omega f(x)\ln f(x){\rm d}\mu(x).$$
With this sign convention for the entropy, the inequalities we derive are of superadditive type; however, the terminology
``subadditivity of the entropy'' is too well entrenched to use anything else.

Now let $p:\Omega \to \R$ be measurable. 
Let $\nu$ be a Borel measure on $\R$, and define $f_{(p)}$ to be the probability density on 
$(\R,{\cal B},\nu)$ such that for all bounded continuous functions $\phi$ on $\R$,
\begin{equation}\label{mardef}
\int_\Omega \phi(p(x))f(x)\dd \mu(x) = \int_\R \phi(t) f_{(p)}(t)\dd \nu(t) .
\end{equation}
In other words, the measure $f_{(p)}\dd \nu$ is the ``push--forward'' of the measure $f\dd \mu$ under $p$:
$$p\#(f\dd \mu) = f_{(p)}\dd \nu\ .$$
The entropy of $S(f_{(p)})$ is defined just as $S(f)$ was,  except with  $(\R,{\cal B},\nu)$ replacing
$(\Omega,{\cal S},\mu)$. We shall be concerned with the following two questions:
\medskip

\noindent{\it (1)} Given $m$ measurable functions $p_1,\dots, p_m$ on $\Omega$,
and $m$ nonnegative numbers $c_1,\dots,c_m$, is there a finite
constant $D$ such that 
\begin{equation}\label{ababs}
\sum_{j=1}^m c_j S(f_{(p_j)}) \le S(f) + D
\end{equation}
for all probability densities $f$ with finite entropy (i.e. satisfying (\ref{finen}))?
\medskip

\noindent{\it (2)} Given $m$ measurable functions $p_1,\dots, p_m$ on $\Omega$,
and $m$ nonnegative numbers $c_1,\dots,c_m$, is there a finite
constant $D$ such that
\begin{equation}\label{ababbl}
\int_\Omega \prod_{j=1}^m f_j(p_j(x)) \dd \mu(x) \le e^D \prod_{j=1}^m \left(\int_\R f_j^{1/c_j}(t)\dd \nu(t)\right)^{c_j}
\end{equation}
for any $m$ nonnegative functions $f_1,\dots, f_m$ on $\R$?

\medskip

For example, consider the case that $\Omega = \R^n$ with its standard Euclidean structure, and $\mu$ is Lebesgue measure on $\R^n$, while $\nu$ is Lebesgue measure on $\R$.  Suppose that $a_1,\dots,a_m$ are $m$ vectors
that span $\R^n$, and define
$$p_j(x) = a_j\cdot x\ .$$
In this case, if we let $X$ denote a random vector with values in $\R^n$ whose law has the density $f$, then $f_{(p_j)}$ is simply the density of the law of $a_j\cdot X$. If we define the entropy of a random variable to be the entropy of its density, provided it has one, we can rewrite
this Euclidean version of (\ref{ababs}) as
$$\sum_{j=1}^mc_jS(a_j\cdot X) \le S(X) + D  .$$
In case $m=n$, $c_j= 1$ for all $j$, and $\{a_1,\dots,a_n\}$ is an orthonormal basis of $\R^n$, 
then this inequality holds with $D=0$, and is the classical subadditivity of the entropy inequality. 

It is even easier to recognize (\ref{ababbl}) as a classical result in this setting: It becomes
$$\int_{\R^n} \prod_{j=1}^m f_j(a_j\cdot x) \dd \mu(x) \le e^D \prod_{j=1}^m \left(\int_\R f_j^{1/c_j}(t)\dd t\right)^{c_j}\ ,$$
which is the classical Brascamp--Lieb inequality.  A celebrated theorem of Brascamp and Lieb~\cite{BL} says that the best constant $e^D$  in this inequality can be computed by using only centered Gaussian functions as trial functions. A new proof based on optimal mass transport was given by Barthe~\cite{B2} who also gave a characterization (depending on the vectors $a_j$ and the constants $c_j$) of when the constant is finite together with a description of the optimizers in some situations. Carlen, Lieb and Loss~\cite{CLL1} introduced a new approach to the Brascamp-Lieb inequalities based on heat flow (see also~\cite{BC}). These authors also completed the gaps left by Barthe in the description of the optimizers. Bennett, Carbery, Christ and Tao~\cite{BCCT1} used a similar approach to deal with the multidimensional versions of the Brascamp-Lieb inequality (see also~\cite{BCCT2} for a direct approach of the finiteness of the constant $e^D$). The paper~\cite{CLL1} (and~\cite{BCCT1} in the multidimensional setting) develops a ``splitting procedure" that will prove useful in our situation too. But we shall see that working with entropy clarifies many technical points.

There are a number of other examples, besides the classical one in the Euclidean setting, where choices of $(\Omega,\mu)$ and $\nu$ lead to
inequalities of interest. For a second  example,  take $\Omega = S^{n-1}$, the unit sphere in 
$\R^n$, and let $\mu$ be the uniform probability measure on   $S^{n-1}$. Then take
$\nu$ to be the probability measure on $\R$ that is  the law of $u\cdot x$, where $u$ is any unit vector
in $\R^n$, so that for all continuous functions $\phi$,
$$\int_{S^{n-1}}\phi(u\cdot x)\dd \mu(x) = \int_\R \phi(t)\dd \nu(t) .$$
(By the rotational invariance of $\mu$, this does not depend on the choice of $u$.)
Now let $\{e_1,\dots,e_n\}$ denote the standard orthonormal basis in $\R^n$, and define
$p_j(x)$ on $S^{n-1}$ by $p_j(x) = e_j\cdot x$.  Then one has the optimal inequalities
\begin{equation}\label{sphere1}
\sum_{j=1}^n \frac{1}{2}S(f_{(p_j)}) \le S(f)\, ,
\end{equation}
for any probability density $f$ on $(\Omega, \mu)$ with finite entropy, and
\begin{equation}\label{sphere2}
\int_{S^{n-1}}\prod_{j=1}^n f_j(e_j\cdot x)\dd \mu(x) \le \prod_{j=1}^n\left(
\int_{S^{n-1}}f_j^2(e_j\cdot x)\dd \mu(x)\right)^{1/2}  =
\prod_{j=1}^n\left(
\int_{[-1,1]}f_j^2(t)\dd \nu(t)\right)^{1/2}\ ,
\end{equation}
for any $n$ nonnegative functions $f_1,\dots, f_n$ on $[-1,1]$. 
See \cite{CLL1} for the original proofs of (\ref{sphere1}) and (\ref{sphere2}), in which 
(\ref{sphere1}) was deduced from (\ref{sphere2}). See
\cite{BCM} for a different and direct   proof of (\ref{sphere1}).  

Since we are concerned in this paper with the relation between subadditivity of entropy and Brascamp--Lieb type inequalities, it is worth recalling the short argument from   \cite{CLL1} that provided the passage from  (\ref{sphere2}) to (\ref{sphere1}):
Let $f$ be any probability density on ${S^{n-1}}$, and let $f_{(p_1)},f_{(p_2)},\dots,f_{(p_n)}$ be its $n$ marginals, as above. 
Then define another probability density $g$ on ${S^{n-1}}$ by
$$g(x) := \frac{1}{C}
\prod_{j=1}^n f_{(p_j)}^{1/2}(e_j\cdot x)\qquad{\rm where}\qquad C := \int_{S^{n-1}}\prod_{j=1}^n f_{(p_j)}^{1/2}(e_j\cdot x)\dd \mu(x)\ .$$
Then by positivity of the relative entropy (Jensen's inequality), we have
\begin{eqnarray}\label{prin1}
0 &\le& \int_{S^{n-1}}\ln\left(\frac{f}{g}\right)f\dd \mu =  S(f) -   \int_{S^{n-1}}\left(\sum_{j=1}^n \ln f_{(p_j)}^{1/2}(e_j\cdot x)\right)f(x)\dd \mu(x) + \ln C\nonumber\\
&=& S(f) -   \frac{1}{2}\int_{\R}\left(\sum_{j=1}^n f_{(p_j)}  \ln f_{(p_j)}\right)\dd \nu  + \ln C\nonumber\\
&=& S(f) -   \frac{1}{2}\sum_{j=1}^n S( f_{(p_j)}) + \ln C\, . 
\end{eqnarray}
Finally, (\ref{sphere2}) implies that
$$C =   \int_{S^{n-1}}\prod_{j=1}^n f_{(p_j)}^{1/2}(e_j\cdot x)\dd \mu(x) \le 
\prod_{j=1}^n\left(\int_{S^{n-1}}f_{(p_j)}(e_j\cdot x)\dd \mu(x)\right)^{1/2} =1$$
since each $f_{(p_j)}$ is a probability density. Thus, $\ln(C) \le 0$, so that (\ref{sphere2}) now follows from (\ref{prin1}).
This argument may give the impression that  (\ref{sphere2}) is a ``stronger'' inequality than (\ref{sphere1}), but as we shall see, this is not the case.

For a third  example,  take $\Omega = {\cal S}_n$, the symmetric group on $n$ letters. 
Let $\mu$ be the uniform probability measure on   $ {\cal S}_n$, and take
$\nu$ to be the uniform probability measure on $\{1,2,\dots,n\}$, so $\nu(i) = 1/n$ for all $i$. 
Define the functions $p_j:\Omega \to  \{1,2,\dots,n\} \subset \R$ by
$p_j(\sigma) = \sigma(j)$ 
for any permutation $\sigma$ of  $\{1,2,\dots,n\}$
Then one has the optimal inequalities
\begin{equation}\label{sphere2b}
\sum_{j=1}^n \frac{1}{2}S(f_{(p_j)}) \le S(f)\, ,
\end{equation}
for any probability density $f$ on $(\Omega, \mu)$, and
\begin{equation}\label{sphere3c}
\int_{{\cal S}_n}\prod_{j=1}^n f_j(p_j(\sigma))\dd \mu(\sigma) \le \prod_{j=1}^n\left(
\int_{{\cal S}_n}f_j^2(p_j(\sigma))\dd \mu(\sigma))\right)^{1/2}
=\prod_{j=1}^n\left(
\sum_{i=1}^nf_j^2(i) \nu(i)\right)^{1/2} \ ,
\end{equation}
for any $n$ nonnegative functions $f_1,\dots, f_n$ on $\{1,\ldots, n\}$. 
See \cite{CLL2} for the  proof of (\ref{sphere3c}).   One could then derive  (\ref{sphere2b}) using the exact same argument that
was used to derive  (\ref{sphere1}) from (\ref{sphere2}).

There are more
examples of interesting specializations of (\ref{ababs}) and (\ref{ababbl}). However, these examples suffice to
illustrate the  context in which the present work is set, and we now turn to the results. 
One basic result of this paper is the following:

\medskip
{\it 
The two questions concerning (\ref{ababs}) and (\ref{ababbl}) that were raised above
are in fact one and the same: We shall prove here that the answer to one question is ``yes'' if and only if the answer to the other question is ``yes'' --- with the same constant
$D$, and with a complete correspondence of cases of equality.}
\medskip

Thus, if one's goal is to
prove a generalized Brascamp-Lieb type inequality, one possible route is to directly prove
the corresponding  generalized  subadditivity of the entropy inequality.  We shall demonstrate the utility of this approach by giving a  simple proof of the classical Brascamp-Lieb inequality on $\R^n$, including a determination of all
of the cases of equality,  through a direct analysis of the entropy. We shall use rather elementary properties of the entropy (scaling properties and conditional entropy) together with geometric properties of the Fisher information.
Moreover, the 
generalized subadditivity of the  entropy inequality that we prove here is new (in its full generality), and is interesting in and of itself.  As we shall see, it turns out to have a rich geometric structure. 
From the point of view of information theory, it might   also be of interest to use the 
 converse implication and to reinterpret  some Brascamp-Lieb inequalities (such as the sharp
Young's convolution inequality) in terms of subadditivity inequalities for the entropy. 

The rest of the paper is organized as follows. In  Section  2, we give the proof that (\ref{ababs}) and (\ref{ababbl}) are dual to one another, so that once one has one inequality established with the cases of equality determined, one has the same for the other. We shall state this duality in a very general setting.

In Section 3, we prove the sharp version of
the general Euclidean subadditivity of the entropy inequality.  

In Section 4 we shall deduce some interesting consequences from this, including a generalization of Hadamard's inequality for the determinant.  

The final Section 5 gives another duality result showing that 
the superadditivity inequalities for Fisher information are dual to certain convolution type inequalities 
of ground state eigenvalues of Schr\"odinger operators.
%
%
These inequalities appear to be new. They may be of some intrinsic interest, but our interest in them here is that a direct proof of the eigenvalue inequalities would yield a direct proof of
Fisher information inequalities that would in turn yield entropy and Brascamp-Lieb inequalities.


\section{Duality of the Brascamp--Lieb inequality and subadditivity of the entropy}

We show that the Brascamp--Lieb inequality is dual to the subadditivity of the entropy, so that once one has proved one 
of these inequalities with sharp constants, one has the other with sharp constants too. In fact, we shall see that there is an
exact correspondence also for cases of equality, but in the next theorem, we focus on the constants.

We shall state the result in a more general setting than the one described in the introduction.
We consider a reference measure space $(\Omega, \mathcal S , \mu)$ and  a family of measure spaces $(M_j, \mathcal M_j , \nu_j)$ together with measurable functions $p_j:\Omega \to M_j$, $j\le m$. For a probability density $f$ on $\Omega$ (with respect to $\mu$), the marginal $f_{(p_j)}$ is thus defined as the probability density on $M_j$ (with respect to $\nu_j$) such that
\begin{equation}\label{mardef2}
\int_\Omega \phi(p_j(x))f(x)\dd \mu(x) = \int_{M_j} \phi(t) f_{(p_j)}(t)\dd \nu_j(t)\, .
\end{equation}
for all bounded measurable functions $\phi$ on $M_j$ ; accordingly the entropies are given by
$$S(f)= \int_\Omega f\ln(f) \, d\mu \quad \textrm{and}\quad S(f_{(p_j)}) = \int_{M_j} f_{(p_j)} \ln(f_{(p_j)}) \dd \nu_j \, .$$
As explained in the introduction, we are mainly interested in the case $(M_j, \mathcal M_j , \nu_j) = (\R, \mathcal B , \nu)$ for all $j\le m$, where $\nu$ is the Lebesgue measure on $\R$.

\medskip

\begin{thm}\label{dual}    
Let $(\Omega,{\cal S},\mu)$ be a measure space, $m\ge 1$  and for $j\le m$, let $(M_j, \mathcal M_j , \nu_j)$ be a measure space together with  a measurable function $p_j$ from $\Omega$ to $M_j$.  For any probability density $f$
on $\Omega$, let $f_{(p_j)}$ the probability density on $M_j$ be defined as in (\ref{mardef2}). Finally, let $\{c_1,\dots,c_m\}$
be any set of $m$ nonnegative numbers.

Then for any $D\in \R$, the following two assertions are equivalent:
\begin{enumerate}
\item For any $m$ nonnegative functions $f_j : M_j \to \R_+$, $j\le m$, we have 
\begin{equation}\label{abBL}
\int_\Omega \, \prod_{j=1}^m f_j(p_j(x)) \dd \mu(x) \le e^D \prod_{j=1}^m \left(\int_{M_j} f_j^{1/c_j}(t)\dd \nu_j(t)\right)^{c_j} \, .
 \end{equation}
\item For every probability density $f$ on  $(\Omega,{\cal S},\mu)$ with finite entropy, we have
\begin{equation}\label{absubH}
\sum_{j=1}^m c_j \, S(f_{(p_j)}) \le S(f) + D \, .
\end{equation}
\end{enumerate}
\end{thm}

The proof  depends an a well known expression for the entropy as a Legendre transform: For any probability density $f$ in $\Omega$, and any function $\phi$ such that
$e^{\phi}$ is integrable,

$$\int_{\Omega} f\ln\left(\frac{e^{\phi}}{f}\right)\dd \mu = \int_{\Omega}f\phi \dd \mu - \int_{\Omega} f\ln f \dd \mu\, .$$
On the other hand, by Jensen's inequality,
$$ \ln\left(\int_{\Omega}e^{\phi}\dd \mu\right)\ge \int_{\Omega} f\ln\left(\frac{e^{\phi}}{f}\right) \dd \mu\, .$$
Therefore,
\begin{equation}\label{ful}
  \int_{\Omega} f\ln f \dd \mu +   \ln\left(\int_{\Omega}e^{\phi}\dd \mu\right) \ge \int_{\Omega}f\phi \dd \mu\, ,
 \end{equation}
and there is equality if and only if $e^{\phi}$ is a constant multiple of $f$ on the support of $f$. We shall use that this Legendre duality nicely combines with the operation of taking marginals.

\medskip

\noindent{\bf Proof of Theorem \ref{dual}:}   First, assume (\ref{abBL}). Consider
 any probability density 
$f$ on $\Omega$, and any $m$ functions $\phi_j$ on  $M_j$, $j\le m$. Using~\eqref{ful} with $\phi$ defined on   $\Omega$ by
\begin{equation}\label{phidef}
\phi(x) := \sum_{j=1}^m c_j\phi_j(p_j(x))\, 
\end{equation}
and~(\ref{mardef2}) we get 
\begin{eqnarray}\label{for}
\int_{\Omega}f(x) \ln f (x)\dd \mu &\ge& \int_{\Omega} f(x)\left(\sum_{j=1}^m c_j\phi_j(p_j(x))\right)\dd \mu  - \ln\left(\int_{\Omega} \prod_{j=1}^m e^{c_j\phi_j(p_j(x))} \dd \mu(x)\right)\nonumber\\
&=& \sum_{j=1}^m c_j\int_{M_j} f_{(p_j)}(t)\phi_j(t)\dd \nu_j(t) - \ln\left(\int_{\Omega} \prod_{j=1}^me^{c_j\phi_j(p_j(x))} \dd \mu(x)\right)\, . \nonumber\\
\end{eqnarray}
Then from the assumption~(\ref{abBL}) applied with $f_j=e^{\phi_j}$, 
$$\int_{\Omega} \prod_{j=1}^me^{c_j\phi_j(p_j(x))} \dd \mu(x)  \le e^D\prod_{j=1}^n\left( \int_{M_j} e^{\phi_j(t)} \dd \nu_j(t) \right)^{c_j}\ .$$
Therefore, (\ref{for}) becomes
\begin{equation}\label{key3}
\int_{\Omega}f(x) \ln f (x)\dd \mu(x) \ge  \sum_{j=1}^m c_j\left( \int_{M_j} f_{(p_j)}(t)\phi_j(t)\dd \nu_j(t) - \ln\left( \int_{M_j} e^{\phi_j(t)} \dd \nu_j(t)\right)\right) -  D\ .
\end{equation}
Now the optimal choice $\phi_j = \ln f_{(p_j)}$ leads to (\ref{absubH}).

Conversely, suppose that  (\ref{absubH}) is true. Consider
$m$ functions $\phi_j$ on $M_j$, $j\le m$, and define $\phi$ on $\Omega$ as in (\ref{phidef}).  Suppose that $e^{\phi}$ is integrable, and  choose $f$ to be
the probability density
\begin{equation}\label{fspec}
f(x) = \left(\int_{\Omega}e^{\phi(x)}\dd \mu(x)\right)^{-1} e^{\phi(x)}\ ,
\end{equation}
so that there is equality in (\ref{ful}). Then we have   from (\ref{ful}) that 
\begin{eqnarray}\label{back}
 \ln\left(\int_{\Omega} \prod_{j=1}^m e^{c_j\phi_j(p_j(x))} \dd \mu(x)\right) &=& 
  \int_{\Omega} f(x)\left(\sum_{j=1}^nc_j\phi_j(p_j(x))\right)\dd \mu(x)  -\int_{\Omega}f(x) \ln f (x)\dd \mu(x)\nonumber\\
&=& \sum_{j=1}^m c_j\int_{M_j} f_{(p_j)}(t)\phi_j(t)\dd \nu_j(t)  -\int_{\Omega}f(x) \ln f (x)\dd \mu(x)\nonumber\\
\end{eqnarray}
On the other hand, (\ref{absubH}) reads as 
\begin{equation}\label{key1}
\int_{\Omega}f(x) \ln f (x)\dd \mu(x) \ge \sum_{j=1}^n c_j \int_{M_j} f_{(p_j)}(t)\ln f_{(p_j)}(t)\dd \nu_j(t) -  D\ ,
\end{equation}
and so (\ref{back}), and then (\ref{ful}) applied on $(M_j,\nu_j)$ with the probability density $f_{(p_j)}$ and the function $\phi_j$ for each $j\le m$, imply
\begin{eqnarray}\label{back2}
\ln\left(\int_{\Omega} \prod_{j=1}^m e^{c_j\phi_j(p_j(x))} \dd \mu\right) &\le& 
\sum_{j=1}^m c_j\left(\int_{M_j} f_{(p_j)}(t)\phi_j(t)\dd \nu_j(t)  - \int_{M_j} f_{(p_j)}(t)\ln f_{(p_j)}(t)\dd \nu_j(t)\right) + D\nonumber\\
&\le& 
\sum_{j=1}^m c_j \ln \left(\int_{M_j} e^{\phi_j(t)}\dd \nu_j(t)  \right) +  D\ .\nonumber\\
\end{eqnarray}
Exponentiating both sides, we obtain (\ref{abBL}). \qed
\medskip

We next examine the relation between cases of equality in the two inequalities.

\begin{thm}\label{cases} Using the
notation  of the previous theorem, 
suppose that $f$ is a probability density on $\Omega$ for which equality holds in the subadditivity inequality (\ref{absubH}).
Then the marginals  $f_{(p_1)}, f_{(p_2)},\dots,f_{(p_m)}$ of $f$  yield
equality in the Brascamp--Lieb  inequality  (\ref{abBL}), and moreover, $f$ and its marginals satisfy
\begin{equation}\label{extr1}f =e^{-D}\prod_{j=1}^m( f_{(p_j)}(p_j(x)))^{c_j}\ .
\end{equation}

Conversely, suppose that  $f_1,\dots,f_m$ are $m$ probability densities (on $M_j$ with respect to  $\nu_j$ for $j=1,\ldots, m$, respectively) for which  equality holds in the Brascamp--Lieb inequality (\ref{abBL}). Then the probability density $f$ defined on $\Omega$ by
$$f(x) := e^{-D} \prod_{j=1}^n( f_j(p_j(x)))^{c_j} $$
 yields equality in the subadditivity inequality  (\ref{absubH})
and moreover $f_j$ is the $j$th marginal of $f$; i.e.  
$f_j = f_{(p_j)}$  for $j\le m$ .
\end{thm}

\medskip
\noindent{\bf Proof:} Suppose that for some probability density $f$, $\sum_{i=1}^m c_i \, S(f_{(p_i)}) - S(f) =  D$. 
Then with this $f$, we must have equality in the first inequality in (\ref{for}), which comes from (\ref{ful}).
By what we have said about the cases of equality in (\ref{ful}), this means that $\phi$, defined in (\ref{phidef}) is a constant multiple of $\ln f$. 
Moreover, to get equality in (\ref{key3}), we were forced to choose $\phi_j = \ln (f_{(p_j)})$. 
This ensures that (\ref{extr1}) is true. 

Furthermore,  to get equality in our intermediate application of the Brascamp--Lieb inequality, we must have that  $\{f_{(p_1)},\dots, f_{(p_n)}\}$ is a set of extremals for the Brascamp--Lieb inequality.  

The other assertion follows in the same way.  \qed

\medskip

By what we have just established, one could try to prove the classical Brascamp--Lieb inequality by first proving a general subadditivity of the entropy inequality for random variables in $\R^n$. We do this in the next section, and shall see that the determination of all of the cases of equality is particularly transparent via this route.  While the Brascamp--Lieb inequality
and subadditivity inequality are equivalent, there is an extra richness 
to the investigation of the cases of equality in
the subadditivity inequality, as this involves statistical independence in a crucial way. Some hint
of this can be seen in the following simple example, which sets the stage for the next section:

Let $m=n$, $c_j= 1$ for all $j$, and $\{a_1,\dots,a_n\}$ be  an orthonormal basis of $\R^n$. 
Take all reference measures to be Lebesgue measure. Then  the Brascamp-Lieb inequality reduces to  an equality, by Fubini's theorem, with $D=0$, and \emph{any}  set of  non negative integrable functions $\{f_1, \ldots, f_n\}$ provides a case of equality. 

On the other hand the dual inequality,  is the \emph{classical subadditivity} of the entropy inequality
$$\sum_{i=1}^m S(X\cdot a_i) \le S(X)\ ,$$
and equality occurs \emph{exactly} when the coordinates $\{X\cdot a_1, \ldots, X\cdot a_n\}$ form a set of \emph{independent} random variables.

In this example, it may appear that the entropy inequality is the more complicated of the two
inequalities. However, the fact that statistical independence enters the picture on the entropy side
is quite helpful: We will make much use of simple entropy inequalities that are saturated only
for independent random variables in our investigation of the cases of equality in the next section.


\section{The general subadditivity of the entropy inequality in $\R^n$}

Let  $\R^n$ be equipped with its standard Euclidean structure. Let $X$ denote a random vector (or variable if $n=1$) with values in $\R^n$, and suppose that $X$ has a density
$f$. 
We denote this correspondence between the random variable $X$ and its density $f$ by writing
$X \sim f$ and set
$$S(X)=S(f)= \int_{\R^n} f(x)\ln f(x) \dd^n x\, .  $$

Thus, in this section, we are specializing the general context of the introduction to the case in which $\Omega$ is $\R^n$, and $\mu$ is Lebesgue measure. We shall also take $\nu$
to be Lebesgue measure on $\R$.

Given a non zero vector $a$ on $\R^n$, identify $a$ with the linear functional 
$a(x) = a\cdot x$.   Then, if $f\sim X$ is a probability density on $\R^n$, $f_{(a)}$, as defined 
by (\ref{mardef}), is the density of $a\cdot X$, that is $f_{(a)}\sim a\cdot X$, and
$$S(X\cdot a)=S(f_{(a)})= \int_{\R} f_{(a)}(t)\ln f_{(a)}(t)\dd t \, .  $$
Note that (\ref{mardef})
specializes to the requirement that 
for every bounded and continuous $\phi:\R\to\R$,
\begin{equation}\label{def:marginal}
\int_{\R^n} \phi (x\cdot a) f(x) \, \dd^n x = \int_{\R} \phi (t) f_{(a)}(t) \, \dd t \ .
\end{equation}
It follows that for all $t\in \R$,
$ f_{(a)}(t) = \frac{1}{|a|}\int_{\{a\cdot x = t\}} f (x)\dd^{n-1} x$.
It is a direct consequence of (\ref{def:marginal}) that  for all $\lambda > 0$, 
\begin{equation}\label{scaling}
f_{(\lambda a)} (t)=\lambda^{-1}f_{(a)}(\lambda^{-1}t)\ .
\end{equation}

With these preliminaries out of the way, we turn to the main question to be addressed in this section: Consider $m$ non zero vectors $a_1,\dots,a_m$ in $\R^n$,
 and $m$ numbers $c_1,\dots,c_m$ with $c_j>0$ for all $j$. 
Then,  we ask:
\medskip

 {\it Is there a finite constant $D\in \R$ so that
 \begin{equation}\label{subadd}
 \sum_{j=1}^m c_j S(a_j\cdot X) \le S(X) +  D
 \end{equation}
 for all random vectors $X$ in $\R^n$, and if so, what is the least such value of $D$, and what are the cases of equality?}
 \medskip

In general there is no finite constant $D$ for which (\ref{subadd}) is true for all $X$. 
There are some simple requirements on $\{a_1,\dots,a_m\}$ 
and $\{c_1,\dots,c_m\}$ for this to be the case. 

 First of all,  for (\ref{subadd}) to hold for any finite constant $D$, the set of vectors
 $\{a_1,\dots,a_m\}$ must span $\R^n$. The following construction is useful for  this and other purposes:   Let $V$ be any proper subspace of $\R^n$, and let $V^\perp$ be its orthogonal complement.
 Then for any number $\lambda>0$, let
 $X_{V,\lambda}$ denote the centered Gaussian random vector (see below for definition) such that 
 \begin{equation}\label{def:X_V}
 \forall u\in V, \ \mathbb E \big((u\cdot X_{V,\lambda})^2\big) = \lambda \quad \textrm{ and }\quad 
  \forall u \in V^\perp ,\ 
 \mathbb E \big((u\cdot X_{V,\lambda})^2\big) = 1.
 \end{equation}
 Then
 \begin{equation}\label{split1}
 S(X_{V,\lambda}) = -\frac{n}{2}\ln(2\pi e) -\frac{{\rm dim}(V)}{2}\ln(\lambda)
 \end{equation}
 while for any $a$ in $\R^n$,
\begin{equation}\label{split2}
S(a\cdot X_{V,\lambda}) = -\frac{1}{2}\ln(2\pi e)  - \frac{1}{2}\ln(\lambda|Pa|^2 + |P^\perp a|^2)\ ,
\end{equation}
where $P$ is the orthogonal projection onto $V$, and $P^\perp = I-P$.

Now take $V$ to be the orthogonal complement of the span of  $\{a_1,\dots,a_m\}$. If the latter is a proper subspace of 
$\R^n$,  then ${\rm dim}(V) \ge 1$, and we see that for any finite $D$,   (\ref{subadd}) would be violated for sufficiently large $\lambda$, since then $|Pa_j|^2 =0$ for each $j$.

 Beyond this spanning condition, there are some simple compatibility conditions that must be
 satisfied by the vectors $a_j$ and the numbers $c_j$. First of all, it follows from (\ref{scaling}) that
 for all $\lambda > 0$,
 $$S(\lambda X) = S(X) - n\ln(\lambda)\qquad{\rm and}\qquad 
 S(a\cdot \lambda X) = S(a\cdot X) - \ln(\lambda)\ .$$
 Therefore, (\ref{subadd}) can only hold when
 \begin{equation}\label{sumcon}
 \sum_{j=1}^mc_j = n\ .
 \end{equation}
 
 There is a further necessary condition  that is somewhat less obvious. The key observation to make is that the right hand side of
 (\ref{split2}) tends to infinity as $\lambda$ tends to zero if and only if $|P^\perp a|^2= 0$, 
 
 Consider any subset
 $J$ of $\{1,\dots,m\}$,  and let 
$$V_J := \textrm{span}\{a_j \; ; \ j\in J\}.$$
Let  $G_J$ denote the Gaussian random variable $X_{V_J,\lambda}$ defined by~\eqref{def:X_V} when $V=V_J$. 
Note that  for each $j\in J$, $|P^\perp a_j|^2 = 0$,  so that for such $j$,
$$S(a_j\cdot G_J) = -\frac{1}{2}\ln(2\pi e) - \frac{1}{2}\ln(|a_j|^2) - \frac{1}{2}\ln(\lambda)\ ,$$
which tends to infinity as $\lambda$ tends to zero. 
Therefore, letting $\lambda$ approach zero, we see that  the leading term in 
$\sum_{j=1}^m c_j S(a_j\cdot G_J)  - S(G_J)$ is at least
$$\frac{1}{2}\left({\rm dim}(V_J) - \sum_{j\in J}c_j\right)\ln(\lambda)\ .$$
(It is exactly this unless for some $i\notin J$, $a_i\in V_J$, in which case we could have taken an even ``worse'' set $J$.)
Hence, if ${\rm dim}(V_J) - \sum_{j\in J}c_j< 0$, there can be  no upper bound on
$\sum_{j=1}^m c_j S(a_j\cdot G)  - S(G)$.
Therefore, 
(\ref{subadd}) can only hold when it is the case that  for all $J$,
 \begin{equation}\label{cond}
 \sum_{j\in J}c_j \le {\rm dim}(V_J)\ .
 \end{equation}
In particular, we must have $c_j \le 1$ for all $j$. 
 
We shall give a simple proof that these necessary conditions are sufficient.  
The following notation shall be used throughout the proof: Given any family 
$\{a_1,\ldots a_m\}$ of vectors spanning $\R^n$, let 
$$A = [a_1,\dots,a_m]$$
denote $n\times m$ matrix whose
$j$th column is $a_j$.  We shall also use $A$ to denote the family 
$\{a_1,\ldots a_m\}$ of spanning vectors.
Thinking of $A$ as the matrix of a linear transformation, 
computed in the canonical bases of $\R^n$ and $\R^m$, will be useful in the proofs of several lemmas below. Note that $A$ has full (row) rank.  Next, let $c$ denote the vector in $\R^m$
 whose $j$th entry is $c_j$.  Finally, define the quantity $D(A,c)$ by
 \begin{equation}\label{ddef}
 D(A,c) := \sup_X\left\{\   \sum_{j=1}^m c_j S(a_j\cdot X) - S(X) \ \right\}\, , 
 \end{equation}
where the supremum is taken over all random vectors $X$ with values in $\R^n$ and with finite entropy. A random vector $X$ for which this supremum is attained will be said to be~\emph{extremal} and will be called an~\emph{extremizer}.

Notice that with $A$ fixed, $D(A,\cdot)$ is the pointwise supremum of a set of affine functions, and as such, it is convex.   We introduce
\begin{equation}\label{def:K_A}
K_A :=\Big\{ c\in [0,1]^m \; ; \ c \textrm{ verifies~\eqref{sumcon} and~\eqref{cond} } \forall J\subset \{1,\ldots, m\} \Big\} ,
\end{equation}
which is clearly a convex subset of the hyperplane of $\R^m$ defined by (\ref{sumcon}). As we have seen, 
$D(A,c)$ is infinite outside $K_A$.  We shall also need later to distinguish the interior of $K_A$ relative to the 
intersection of $[0,1]^m$ and the hyperplane specified by~\eqref{sumcon}:
\begin{equation}\label{def:K_Aint}
K_A^\circ :=\left\{ c\in K_A \; ;  \sum_{j\in J}c_j < {\rm dim}(V_J), \ \forall J\subsetneq \{1,\ldots, m\}, \; J\neq \emptyset \right\}.
\end{equation}
We shall make an extensive use of the fact that $K_A$ and $K_A^\circ$ are invariant under linear transformation, in the sense that for any invertible linear operator $T$ on $\R^n$, we obviously have $K_{TA}=K_A$ and $K_{TA}^\circ = K_{A}^\circ$ 
 with the notation $TA=[Ta_1, \ldots, Ta_m]$ when $A=[a_1, \ldots, a_m]$.

Also define $D_{\cal G}(A,c)$, the Gaussian analog of~\eqref{ddef}, by
 \begin{equation}\label{dgdef}
 D_{\cal G} (A,c) := \sup_G\left\{\   \sum_{j=1}^m c_j S(a_j\cdot G) - S(G) \ \right\}\ .
 \end{equation}
in which the supremum is taken over all centered Gaussian random vector $G$ with values in $\R^n$.  By a \emph{centered Gaussian random vector}, we mean one that has a density of the form
$$\frac{1}{|\det(C)|}\left(\frac{1}{2\pi}\right)^{n/2}e^{-|C^{-1}x|^2/2}$$
for some symmetric invertible matrix $C$ on $\R^n$. More generally, a \emph{Gaussian random vector} is a random vector of the form $x_0+G$ with $x_0\in \R^n$ and $G$ a centered Gaussian random vector. We can restrict ourselves to centered random vectors because the entropy is invariant under translation. A Gaussian random vector is said to be \emph{isotropic} if its covariance matrix is a multiple of the identity; it is said to be \emph{standard} if it is centered and if its covariance matrix is the identity (i.e.  it is a $\mathcal N (0, \id)$ Gaussian vector). 

At this point, it is important to note that all the definitions made so far make sense more generally on a finite dimensional Euclidean space $(E,\cdot)$. We have made the identification $E=\R^n$, which has the advantage to allow us to work with matrices. Later, we shall also need to work on subspaces of $\R^n$, which are then canonically equipped with the Euclidean structure inherited from $\R^n$; we then need to work with the corresponding Euclidean versions of the notions introduced above.

It is clear that $D_{\cal G} (A,c)$ is also a convex function of $c$, and that 
$D_{\cal G} (A,c)\le D(A,c)$. Also, since our proof that $D(A,c) = \infty$
for $c\notin K_A$ used a centered Gaussian random vector, it shows also  that
$D_{\cal G} (A,c) = \infty$
for $c\notin K_A$.  In fact, we have the following:

\medskip
\begin{thm}\label{gaussian} For every family $A=\{a_1, \ldots, a_m\}$ of $m$ vectors spanning $\R^n$ and every vector $c$ in $\R^m$ 
with $0\le c_j\le 1$ for all $j$,  we have
$$D(A,c) = D_{\cal G} (A,c)  \ ,$$
and furthermore $D(A,c)$ is finite if and only if $c\in K_A$. 
\end{thm}
\medskip

The proof will be accomplished in three steps:
\medskip

\noindent{\it Step 1:}  We shall first consider the case in which the vectors $a_j$ are all unit vectors $u_j$ satisfying the following special condition, put forward by K. Ball in the setting of Brascamp-Lieb inequalities (see e.g.~\cite{Ball}): 
  \begin{equation}\label{fj}
  \sum_{j=1}^m c_j\, u_j\otimes u_j = \id_{\R^n}\, ,
  \end{equation}
with $c_j\ge 0$. (Note that~\eqref{sumcon} automatically holds, as it can be seen by taking the trace, and that $c_j\le1$ for all $j\le m$.)
Under this condition, we give a simple proof of Theorem \ref{gaussian} using an elementary superadditivity property of the  Fisher information and integration along the heat flow.
The proof here draws on ideas from \cite{BCM}.

\medskip

\noindent{\it Step 2:}   We shall show that for $c\in K_A^\circ$, there is a linear change of variables
that reduces this case to the one considered in the first step.  While the lemma that provides the existence of the change of variables would appear to be a simple statement about linear algebra, the existence of this change of variables is intimately
connected with the existence of Gaussian optimizers for the subadditivity (and hence the Brascamp--Lieb) inequality.

\medskip

\noindent{\it Step 3:}   We show that on  $K_A\backslash K_A^\circ$, the variational problem in (\ref{ddef}) may be ``split'' into two 
problems of the same type, but each involving only a subsets of the original vectors, and integration over a proper subspace of $\R^n$.
Repeating this splitting operation, one eventually reduces to variational problems of the type considered in the second step.
This step is modeled after a similar splitting argument developed in \cite{CLL1}, but as we shall see, the entropic version has advantages that will help us determine all of the cases of equality.

\medskip

\noindent{\bf Remark:} 
If one is content to prove only that  $D(A,c)$ is finite if and only if $c\in K_A$, there is a very expeditious route: One can easily check the finiteness of  $D(A,c)$ at the extreme points of 
$c\in K_A$ (where, as shown by Barthe, each $c_j$ is either $0$ or $1$). Then the convexity of 
$D(A,c)$ implies finiteness on all of $K_A$, and we know it is infinite outside. Proving the equality 
$D(A,c) = D_{\cal G} (A,c)$ on all of $K_A$ is more subtle: The values of  $D(A,c)$ and $
D_{\cal G} (A,c)$ {\em do} jump as one crosses the boundary of $K_A$, and we see nothing to preclude
$D(A,c)$ from jumping up more than $D_{\cal G} (A,c) ${\em on} the boundary.   Thus, it is
not only for the classification of the cases of equality that we argue as we do in the third step:
we do not know of any quick way to ``pass to the boundary'' of $K_A$ and wrap of the proof of
Theorem~\ref{gaussian} after the second step without developing the splitting argument.

\medskip

We now begin with the first step. Here we shall use a 
simple  superadditivity result for the {\it Fisher information}:
If $X\sim f$ is a random vector with a differentiable density $f$, define the Fisher information of $X$ or of $f$ by
\begin{equation}\label{idef}
I(X)=I(f)= \int_{\R^n} \frac{|\nabla f |^2}f\, .
\end{equation}

This quantity is related to the entropy through  the heat flow  as follows: Let $\Delta$ denote the Laplacian on $\R^n$, and let
$G$ denote a standard Gaussian random vector on $\R^n$ independent of $X$, so that if $f\sim X$, 
$$e^{t\Delta}f  \sim X + \sqrt{t}G\, .$$
Then we have the identity
$$\frac{{\rm d}}{{\rm d}t}S(e^{t\Delta}f) = -I(e^{t\Delta}f)\ ,$$
and in particular, the right hand side is finite for all $t>0$. 

The basic inequality concerning the Fisher information that will yield us our subadditivity result
 is the fact that for any unit vector $u$,
\begin{equation}\label{supadd1}
I(f_{(u)})  = I(u\cdot X) \le   \int_{\R^n} \frac{|u\cdot \nabla f |^2}f\ ,
\end{equation}
with equality if and only if $f$ is the product of $f_{(u)}$ and a probability density $g$ on the orthogonal complement of $u$.
This was proved in \cite{C}; see Theorem 2 there with $p=2$. Let us include here for completeness a different proof taken from~\cite{BCLM} (were more abstract settings are studied). This proof requires more regularity than the one in \cite{C}, but that is fine for our purpose,
as we shall apply the inequality along the heat flow. 

Using the definition of the marginal~\eqref{def:marginal} twice and H\"older's inequality, we have:
\begin{eqnarray*}
I(f_{(u)}) &=& - \int_\R f_{(u)} (\ln f_{(u)})''\dd t = - \int_{\R^n} f(x) (\ln f_{(u)})''(x\cdot u)\dd^n x \\
& = &  \int_{\R^n} \frac{(f_{(u)})'(x\cdot u)\, (u\cdot \nabla f(x))}{f_{(u)}(x\cdot u)}\dd^n x \\
& \le & \sqrt{\int_{\R^n} \frac{[(f_{(u)})' (x\cdot u)]^2 }{(f_{(u)}(x\cdot u))^2}\, f(x) \dd^n x\,
\int_{\R^n} \frac{(u\cdot \nabla f(x) )^2 }{f(x)}\dd^n x }   \\
& = & \sqrt{\int_{\R} \frac{[(f_{(u)})']^2 }{(f_{(u)})^2}\, f_{(u)}\dd t\,
\int_{\R^n} \frac{(u\cdot \nabla f )^2 }{f}\dd^n x}\\
& = & \sqrt{I(f_{(u)}) \int_{\R^n} \frac{(u\cdot \nabla f )^2 }{f}\dd^n x }\, .
\end{eqnarray*}
This proves~\eqref{supadd1}. Equality in~\eqref{supadd1} requires equality in H\"older's inequality and so for some $\lambda \in \R$ we have $(u\cdot \nabla)\log f (x) = \lambda (\log f_{(u)})'(x\cdot u)$ for all $x\in \R^n$; this $\lambda$ has to be $1$ for equality to hold in~\eqref{supadd1}  and therefore $f(x)=f_{(u)}(x\cdot u) h(x - (x\cdot u)u)$ for some probability density $h$ on $u^\perp$.

From~\eqref{supadd1}, we immediately deduce the superadditivity of information. But before stating the result, let us make a definition needed to discuss the cases of equality.

\medskip

\begin{defi}[Reducible spanning set]  Let $\{a_1,\dots,a_m\}$ be any set of $m$ vectors spanning $\R^n$. It is a \emph{reducible
spanning set} in case there are two proper subspaces $V_1$ and $V_2$ of $\R^n$ such that
$\R^n = V_1\oplus V_2$, and such that each $a_j$ belongs to either $V_1$ or to $V_2$. Otherwise, 
$\{a_1,\dots,a_m\}$ is called an \emph{irreducible spanning set}.
\end{defi}

\begin{prop}\label{fisher} \label{special}
Consider any set $\{u_1,\dots,u_m\}$ of $m$  \emph{unit} vectors in $\R^n$, such that there are numbers
$\{c_1,\dots,c_m\}$, with $0\le c_j \le 1$ for each $j\le m$,  so that the decomposition of the identity~(\ref{fj}) is satisfied. 
Let $G$ denote a standard Gaussian random vector.

Then for all random vectors $X$ with finite Fisher information,
 \begin{equation}\label{infosub}
 \sum_{j=1}^m c_j \, I(u_j\cdot X) \le I(X)\, ,
 \end{equation}
with equality if $X=G$, and for all random vectors $X$ with finite entropy
\begin{equation}\label{bcn181}
\sum_{j=1}^m c_j\,  S(u_j\cdot X) - S(X)  \le \sum_{j=1}^m c_j \, S(u_j\cdot G) - S(G) = 0  \, .
\end{equation}

Moreover there is equality in these inequalities if and only if for each $j\le m$,
$u_j\cdot X$ and $X - (u_j\cdot X)u_j$ are independent.
Under the condition that  $n\ge 2$ and that
$\{u_1,\dots,u_m\}$ is an irreducible spanning set, then there is equality in these inequalities if and only if $X$ is  an isotropic Gaussian random vector. 
\end{prop}
\medskip

Note that this proposition in particular implies that $D(U,c)=D_\mathcal{G}(U,c)=0$ when $U=[u_1, \ldots , u_m]$ are unit vectors of $\R^n$ and $c=(c_1,\ldots, c_m)$ nonnegative real numbers satisfying~\eqref{fj}. 

The proof of~\eqref{infosub} and~\eqref{bcn181} is elementary and follows~\cite{BCM}. The determination of the cases of equality requires a bit more work, but it remains quiet direct  (compared to analogous result on the side of the Brascamp-Lieb inequality).
\medskip

\noindent{\bf Proof:} Inequality~(\ref{infosub}) follows immediately from (\ref{supadd1}) 
and condition~(\ref{fj}) rewritten in the form
$$\forall x \in \R^n , \quad \sum_{j=1}^m c_j \, (x\cdot u_j)^2 = |x|^2 .$$ 
Equality for $X=G$ is obvious as $G\cdot u_i$ is a standard Gaussian variable and so the computation boils down to the equality $\sum c_j = n$. (For the same reason the right-hand side of the inequality~\eqref{bcn181} is zero.)

 As we have noted, the Fisher information of $f$
is related to the entropy of $f$ through
${\displaystyle \frac{\dd}{\dd t} S(e^{t\Delta} f) = -I(e^{t\Delta}f)}$.
It is also easy to see (using that $\Delta$ commutes with translations) that if $u$ is any unit vector, then $f_{(u)}$, the marginal of $f$ along $u$, has the property that
$(e^{t\Delta}f)_{(u)} = e^{t\Delta}f_{(u)}$ where we keep the same notation of the $1$-dimensional heat semi-group ($\Delta g = g''$ in dimension $1$); we again have  (in dimension $1$) that
$$\frac{\dd}{\dd t} S((e^{t\Delta} f)_{(u)}) = -I((e^{t\Delta}f)_{(u)})\, .$$

Then since $e^{t\Delta}f \sim  X + \sqrt{t}G$, and because  $\sum_{j=1}^m c_j S(u_j\cdot X) - S(X)$ is invariant under dilation, i.e. under the substitution $X \to \lambda X$, we get 
$$\left[\sum_{j=1}^m c_j\, S(u_j\cdot G) - S(G)\right] - 
\left[\sum_{j=1}^m c_j\, S(u_j\cdot X) - S(X)\right] = \int_0^\infty \left[I(e^{t\Delta} f) - \sum_{j=1}^m c_j\,  I((e^{t\Delta} f)_{(u_j)})\right] \dd t\, .$$
By Theorem \ref{fisher}, the integrand above is non negative for all $t$, and so  (\ref{bcn181})
is proved. 

The condition for cases of equality in~\eqref{supadd1} tell us that there is equality in~\eqref{infosub} for a random vector $X$ with finite Fisher information if and only if $X$ verifies the following property $(\mathcal P)$:
$$(\mathcal P) \qquad \forall i\le m, \quad X\cdot u_i \textrm{ and } X- (X\cdot u_i)u_i \textrm{ are independent.}$$
If $G$ is a standard Gaussian random vector independent of $X$, then  $X$ verifies $(\mathcal P)$ if and only if for all $t>0$, $X + \sqrt t G $ verifies $(\mathcal P)$. Thus for a random vector with finite entropy, 
there is equality in~\eqref{bcn181} if and only if $X$ verifies $(\mathcal P)$. 

Our goal is now to characterize, when $n\ge 2$, the random vectors verifying $(\mathcal P)$ under the assumption that $\{u_1, \ldots , u_m\}$ is an irreducible spanning set of unit vectors. First note that if we prove that $X+\sqrt{t}G$ is an isotropic Gaussian  for all $t>0$, then so is $X$. Therefore, using again the stability of the property $(\mathcal P)$, we need only consider random vectors $X$ with smooth and strictly positive density. Secondly, we can assume that no two vectors of the family $\{u_i\}_{i\le m}$ are linearly dependent. Indeed, by keeping only one representative for the subspaces $\R u_j$, we construct a subfamily of the vectors $\{u_i\}_{i\le m}$ which span $\R^n$ and which remains irreducible.

So from now let $\{u_1, \ldots , u_m\}$ is an irreducible spanning set of unit vectors of $\R^n$ ($n\ge 2$), with no two vectors linearly dependent, and $X$ a random vector verifying $(\mathcal P)$ and with a smooth density $f>0$. Thus for every $i\le m$ there exists two probability densities $g_i$ and $h_i$, on $\R$ and $u_i^\perp \simeq \R^{n-1}$ respectively, such that
$$f(x) = g_i(x\cdot u_i) h_i(x- (x\cdot u_i)x)$$
Writing $F=\log f$, $G_i = \log g_i$ and $H_i=\log h_i$ for each $i\le m$, we have
$$F(x) = G_i(x \cdot u_i) + H_i(x- (x\cdot u_i)x) \, ,$$
so that
$$(u_i\cdot\nabla)F(x) = G_i'(u_i\cdot x)\, .$$
Hence for any $j\ne i$,
$$(u_j\cdot\nabla)(u_i\cdot\nabla)F(x)  =  (u_i\cdot u_j) G_i''(u_i\cdot x)\, .$$
Interchanging the roles of $i$ and $j$, 
$$(u_j\cdot\nabla)(u_i\cdot\nabla)F(x)  =  (u_i\cdot u_j) G_j''(u_j\cdot x)\, .$$
Evidently the left hand side depends on $x$ only thorough $u_i\cdot x$ and only through $u_j\cdot x$.
But since $u_i$ and $u_j$ are linearly independent, this means that the left hand side is constant. 
Hence,
\begin{equation}\nonumber
 \text{for\ every\ }\  i\ne j, \  (u_i\cdot\nabla)(u_i\cdot\nabla)F\ \text{ is\  constant.}
 \end{equation} 
Furthermore, under the condition that $\{u_1,\dots,u_m\}$ is an irreducible spanning set, 
 if any one vector $u_i$ is removed from $\{u_1,\dots,u_m\}$, the remaining vectors still span
$\R^n$.  For otherwise, since $m\ge n\ge 2$, we could take $V_1$ to be the span of $\{u_i\}$, and take $V_2$ to be the span of
$\{u_1,\dots,u_m\}\backslash u_i$, and we would have $\R^n = V_1\oplus V_2$. Thus each $u_i$ decomposes in the generating  family $\{u_j\}_{j\neq i}$ and therefore,
\begin{equation}\nonumber
 \text{for\ every\ }\  i, j \le m, \  (u_i\cdot\nabla)(u_i\cdot\nabla)F\ \text{ is\  constant.}
 \end{equation} 
But this implies that the Hessian of $F$ is constant. Thus, $X$ is Gaussian. 
To prove that this Gaussian is isotropic, let $C$ be the covariance matrix of $X$. Then  property $(\mathcal P)$ 
implies that each $u_i$ is an eigenvector of $C$. Since eigenvectors of symmetric matrices
are orthogonal if they have distinct eigenvalues, all of the eigenvalues must be the same unless there is such a ``splitting'' of $\R^n$
into at least two (orthogonal) subspaces that together contain all of the vectors $u_j$.  This would contradict the hypothesis that
$\{u_1,\dots,u_m\}$ is an irreducible spanning set.
\qed

\medskip

The following lemma will facilitate the application of the the statement concerning the cases of equality in Proposition~\ref{special}:

\medskip
\begin{lm}\label{char} Let $A = \{a_1,\dots,a_m\}$ be any family of  $m$ vectors spanning $\R^n$.
If $\{a_1,\dots,a_m\}$ is a reducible spanning set and $D(A,c)$ is finite, then $c \notin K_A^\circ$. 
\end{lm}
\medskip

\noindent{\bf Proof:} Let $\R^n = V_1\oplus V_2$ be a decomposition of $\R^n$ into two proper subspaces 
such that each $a_j$ is contained in one of them or the other.  Let $V$ be the
{\em orthogonal complement} of $V_1$, $\R^n = V_1\stackrel{\perp}{\oplus} V$ and let  $X_{V,\lambda}$ be the Gaussian random variable defined as in (\ref{def:X_V}).   
Then by  (\ref{split1}) and  (\ref{split2}), with $P$ denoting the orthogonal projection onto $V$,
\begin{eqnarray}
\sum_{j=1}^m c_j S(a_j\cdot  X_{V,\lambda}) - S(X_{V,\lambda}) &=& 
-\frac{1}{2}\left(\sum_{j:a_j\in V_1} c_j\ln(|a_j|^2) + \sum_{j:a_j\in V_2}c_j\ln(\lambda |P a_j|^2 + |P^\perp a_j|^2)\right) \nonumber\\ 
&+& \frac{1}{2}{\rm dim}(V)\ln(\lambda)\, ,
\end{eqnarray}
with $Pa_j\neq 0$ for $j\in V_2$, since $Px =0 \Rightarrow x\in V_1$.
Then, using  that  ${\rm dim}(V) =   {\rm dim}(V_2)$, this expression (in $\lambda$)
has the form $$\frac{1}{2}\left({\rm dim}(V_2)- \sum_{j:a_j\in V_2}c_j\right)\ln(\lambda) \ + \ ({\rm terms\  bounded\  in}\ \lambda >1)\ ,$$
which is unbounded for large $\lambda$ unless
 $$\sum_{j:a_j\in V_2}c_j = {\rm dim}(V_2)\ .$$
This must be the case since by hypothesis that $D(A,c) < \infty$. Thus, $c\notin K_A^\circ$ \qed

\medskip

We have now completed the first step. We start the second by showing that the change of variables matrix $R$ does exist for $c\in K_A^\circ$. 
The existence of such a change of variables can be deduced from results of Bennett-Carbery-Christ-Tao~\cite{BCCT1}. However, the flow of logic in their deduction (and in~\cite{CLL1})  runs counter to ours:
They first show that such a change of variables exists whenever there are Gaussian optimizers for the Brascamp--Lieb problem, and then show that Gaussian optimizers exist for
$c\in K_A^\circ$. Here, we need the change of variables at the outset of
our analysis, and hence need a direct proof of this result. We now
provide one, using a geometric result of Barthe. 
\medskip
\begin{lm}\label{fjcond} Let $A=\{a_1,\dots,a_m\}$ be any family of $m$ vectors that span $\R^n$. 
Let $\{c_1,\dots,c_m\}$ be any $m$ numbers verifying $0\le c_j\le 1$ and satisfying~\eqref{cond}. If  $c\in K_A^\circ$, then there exists an invertible symmetric  $n\times n$ matrix $R$ so that
\begin{equation}\label{propR}
\sum_{j=1}^m c_j\left(\frac{Ra_j}{|Ra_j|}\right)\otimes \left(\frac{Ra_j}{|Ra_j|}\right) = \id_{\R^n} \, .
\end{equation}
When $n\ge 2$, there is exactly one such matrix $R$ satisfying the further requirements that $R$ be positive definite, and that ${\rm trace}(R^2) = n$.
On the other hand, for $c\notin K_A$, no such matrix $R$ exists.
\end{lm}
\medskip
\noindent{\bf Remark:} After settling the cases of equality in Theorem
\ref{gaussian} we shall derive  {\em necessary and sufficient conditions} for the existence of such a matrix $R$.   Though the conditions are simple and explicit, it turns out that the matrix $R$ exists if and only if 
the supremum in (\ref{dgdef}) is attained at some centered Gaussian $G$, and our proof that the conditions we give are necessary
and sufficient depends on this.

\medskip
\noindent{\bf Proof:} 
Take any diagonal $m\times m$ matrix $S$ with positive diagonal entries $s_j$, $j\le m$, and define the $n\times n$ matrix $R_S$ by
 $$R_S= ((AS)(AS)^t)^{-1/2}\ .$$
 This makes sense since $(AS)(AS)^t$ is a positive definite $n\times n$ matrix. Notice that  $(R_SAS)(R_SAS)^t = \id_{\R^n}$,
 or, what is the same
 $$\sum_{j=1}^m s_j^2 R_Sa_j\otimes R_Sa_j = \id_{\R^n} .$$
 Therefore,
 $$\sum_{j=1}^m c_j \left(\frac{s_j}{\sqrt{c_j}} R_Sa_j\right)\otimes  \left(\frac{s_j}{\sqrt{c_j}} R_Sa_j\right) = I\ .$$
 We have what we seek if and only if for each $j$, ${\displaystyle \frac{s_j}{\sqrt{c_j}} R_Sa_j}$ is a unit vector, which is the case if and only if for each $j$,
 $ c_j = s_j^2|R_Sa_j|^2$. By the definition of $R_S$, this means
 \begin{equation}\label{el}
 c_j = e_j\cdot [ (AS)^t ((AS)(AS)^t)^{-1}(AS) ] e_j
 \end{equation}
 where $\{e_1,\dots,e_m\}$ denotes the standard orthonormal basis in $\R^m$. Note that  
$e_j\cdot(AS)^t((AS)(AS)^t)^{-1}(AS)e_j$   is also the $j$th diagonal entry of the orthogonal projection in $\R^m$ onto the image of $(AS)^t$. 
 
It has been shown \cite{B2} (see \cite{CLL1} for another proof and a statement in this formulation) that there exists positive numbers $s_1,\dots, s_m$ for which (\ref{el}) is true
 whenever $c\in K_A^\circ$, and that in this case, when $n\ge 2$, 
 the set of numbers is unique up to a common multiple. Thus, for $c\in K_A^\circ$, such an $R$
 exists.
 
As for the uniqueness, note that given any such matrix  $R$, we can change variables, replacing
$X \to R^{-1}X$ and $a_j \to u_j := |Ra_j|^{-1}Ra_j$. Then Proposition \ref{special} may be applied
to deduce that the only extremizers for the new problem are isotropic Gaussians. Undoing the
change of variables, we see that the only extremizers  of the original problem are Gaussians whose
covariance is a multiple of $R^2$.  Thus, under the further condition that
$R$ be positive definite (instead of simply symmetric), and that the trace of $R^2$ is fixed, $R$ is uniquely determined.
 
The same change of variables argument (which is exploited systematically in  Lemma \ref{interior} below) shows, through Proposition 
~\ref{special},  that if such a matrix $R$ exists, then $D(A,c)<\infty$.
 As we have seen, this is impossible when $c\notin K_A$.
 \qed
 
 \medskip
 
 \noindent{\bf Remark:}  The first proof that there exists a solution, essentially unique, to (\ref{el}) whenever $c\in K_A^\circ$
 is due to Barthe \cite{B2}. However, he used a different characterization of $K_A$, and did not mention the condition (\ref{cond}).
 Another proof of this, based directly on (\ref{cond}) was given in \cite{CLL1}, together with a proof that the characterization
 of $K_A$ in Barthe's paper is equivalent to the one based on (\ref{cond}).  
 
\medskip

With the change of variable provided by the previous lemma, we can finish the second step and describe what happens when $c\in K_A^\circ$.

\medskip
\begin{lm}\label{interior} For any family $A=\{a_1, \ldots, a_m\}$ of $m$ vectors spanning $\R^n$, and all vectors $c$ in $K_A^\circ$,
$$D(A,c) = D_{\cal G} (A,c)  \, ,$$
and there exist a Gaussian optimizer. 
Moreover, if   $n\ge 2$, then $\sum_{j=1}^m c_jS(a_j\cdot X) - S(X) = D(A,c)$ if and only if  $X$ is Gaussian and its covariance is a constant multiple of $R^2$ where $R$ is the unique positive definite matrix   verifying~\eqref{propR}  with $\textrm{Tr}(R^2)=n$.
\end{lm}
\medskip

 \noindent{\bf Remark:} The  condition ``$n\ge 2$", which  has already appeared several times,
 is present because in one dimension, the subadditivity problem is trivial, so that Gaussians play no special role. Indeed, assume we are given $c_1, \ldots, c_m\ge 0$ with the condition that $\sum c_j = 1$ and $A=\{a_1, \ldots, a_m\}$ a family of non-zero  real numbers. Then,
 setting
 $$D:= -\sum_{j=1} ^m c_j \log|a_j| $$
 we have, for every random variable $X$ on $\R$ with finite entropy
 $$\sum_{j=1}^m c_j S(a_j X) - S(X) = D.$$
 Therefore $D(A,c)=D$ and every random variable $X$ is an extremizer.
 
\medskip

\noindent{\bf Proof:}  Let $R$ be an invertible symmetric matrix verifying~\eqref{propR} provided by the Lemma~\ref{fjcond}.
Since  for any random vector $X$ with finite entropy, we have
$$S(X\cdot a_j)= S\left(\frac{Ra_j}{|Ra_j|}\cdot R^{-1} X \right)  - \ln(|Ra_j|) \qquad{\rm and}\qquad
S(X)= S(R^{-1}X) - \ln(|\det(R)|)   \, , $$
we obtain
\begin{eqnarray}
\sum_{j=1}^m c_jS(a_j\cdot X) - S(X) &=&  \sum_{j=1}^m c_j S\left(\frac{Ra_j}{|Ra_j|}\cdot R^{-1} X \right) - S(R^{-1}X)\nonumber\\
&-& \sum_{j=1}^m c_j\ln(|Ra_j|) + \ln(|\det(R)|)\, .\nonumber
\end{eqnarray}
Introduce the family of vectors ${\displaystyle u_j:= \frac{Ra_j}{|Ra_j|}}$ for $j\le m$, and set $U=[u_1,\ldots,u_m]$.
The previous equality implies that
 \begin{equation}\label{trans}
 D(A,c)  = D(U,c)  - \sum_{j=1}^m c_j\ln(|Ra_j|) + \ln(|\det(R)|) \, .
 \end{equation}
 Thus we are reduced to studying the problem determining $D(U,c)$ and the extremizers there (noting  that $X$ is an extremizer for $D(A,c)$ if and only if $R^{-1}X$ is extremizer an for $D(U,c)$).
Note also  that since the  vectors $\{u_1,\dots,u_m\}$ are obtained from the  vectors $\{a_1,\dots,a_m\}$ by a non singular linear transformation, they span $\R^n$, and we have $K_U^\circ=K_A^\circ \ni c$.

Since $U=[u_1, \ldots, u_m]$ is a family of \emph{unit} vectors verifying the decomposition of the identity~\eqref{fj}, we can apply Proposition~\ref{fisher} and get that
\begin{equation}\label{reduced}
D(U,c) =   D_{\cal G}(U,c) =0< \infty\, ,
\end{equation}
and every isotropic Gaussian vector is an extremizer.
To prove that all optimizers are Gaussian when $n\ge 2$,  note first that, by Lemma \ref{char}, $c\in K_U^\circ$  implies that  $\{u_1,\dots,u_m\}$ is an irreducible spanning set. Therefore any optimizer of the variational problem defining $D(U,c)$ is an
isotropic Gaussian.
(Then  every optimizer for $D(A,c)$ is Gaussian whose  covariance is a multiple of $R^2$.)
\qed
\medskip

 \noindent{\bf Remark:} 
 Note that the proof above gives also  the following statement:   If there exists an  invertible matrix $R$  verifying~\eqref{propR} then (with no further assumptions on $c$ and $A$) we have  that  $D(A,c)<+\infty$ and that $RG$ is an extremizer for every   standard Gaussian vector $G$.
\medskip

We now turn to the third step.  When $c\in K_A\backslash K_A^\circ$, we will pick a non-empty proper subset  $J$ of $\{1,\dots,m\}$  of least cardinality among subsets for which equality holds in (\ref{cond}).   We shall now show
that the variational problem defining $D(A,c)$ splits into two such problems involving fewer vectors and random variables in a lower
dimensional space. Repeated splittings, and what we have already proved, will enable us to settle all questions concerning
the variational problem   defining $D(A,c)$.
The splitting argument presented here is patterned on one developed in \cite{CLL1} for the Brascamp--Lieb inequality. However, as we shall see, in the subadditivity setting, the argument leads to a clear and simple analysis of cases of equality.  It relies on properties of the conditional entropy. 

As mentioned at the beginning of this section, we shall need to work on subspaces of $\R^n$ and thus make use of the definition made above in the setting of Euclidean spaces. For a given family $A=\{v_1, \ldots v_k\}$ of vectors on $\R^n$, we introduce the Euclidean subspace $E:={\rm span}(v_1, \ldots, v_k)$ equipped with the induced Euclidean structure from $\R^n$ (i.e. the scalar product is the same). For real numbers $c_1, \ldots, c_k$ with $0\le c_j \le 1$, the quantities $D(A,c)$ and $D_G(A,c)$ are then implicitly assumed to be defined on the Euclidean subspace $E$ (the random vectors live on $E$ and the entropies are computed with respect to the Lebesgue measure on $E$, where the laws of the vectors live). Accordingly, the set $K_A$ is to be understood as
$$K_A :=\Big\{ c\in [0,1]^{k} \; ; \  \sum_{j=1}^k c_j = \textrm{dim}(E)  \textrm{ and~\eqref{cond} holds } \forall J\subset \{1,\ldots, k\} \Big\} .$$

Let us fix the following notation.  Let  $A=\{a_1,\ldots, a_m\}$  be a family of of $m\ge 1$ vectors spanning an Euclidean space $E$,
$$E=  {\rm span}\left(\{a_j\ ;\ j\in I\}\right)\, $$
(in a first step we shall have $E=\R^n$). For family of $m$ real numbers $c\in K_A$ and a non-empty proper subset  $J$ of $\{1,\dots,m\}$ for which equality holds in~(\ref{cond}), denote by $P_J$ the orthogonal projection onto $V_J=\textrm{span}\{a_j\; ; \ j\in J\}$ and let $P_J^\perp = \id_E -P_J$ be the complementary projection. 
Define, for $j\in J^c:=\{i\in\{1,\ldots,m\}\; ; i\notin J\}$ the vector $b_j=P_J^\perp a_j$ and 
$$A_J=[a_j \, ; j\in J] \quad \textrm{ and }\quad B_{J^c}=[b_j \, ; j\in J^c]$$
the ordered (by ordering $J$ and $J^c$ as increasing subsequences of $1,\ldots,m$) families of vectors $(a_j)_{j\in J}$ and $(b_j)_{j\in J^c}$.
For any subset $K$ of  $\{1,\dots,m\}$, and $c\in \R^m$, let $c_K$ denote the vector of $\R^{|K|}$ whose coordinate are the $(c_j)_{j\in K}$ ($K$ being written as an increasing subsequence of $1,\ldots, m$).
Since there is equality in~\eqref{cond} for $J$, we have
$$c_J \in K_{A_J}.$$
Note that $V_J + V_{J^c}=E$ ({\it a priori} this sum is not direct) and so $V_{J}^\perp=P_J^\perp V_{J^c}$. Thus we have $V_J^\perp={\rm span}\left(\{b_j\ :\ j\in J^c\}\right) $, i.e.:
\begin{equation}\label{equal2c}
E= {\rm span}\left(\{a_j\ ;\ j\in J\}\right) \stackrel{\perp}{\oplus}  {\rm span}\left(\{b_j\ ;\ j\in J^c\}\right)\, .
\end{equation}
And  we also have
$$c_{J^c}\in K_{B_{J^c}}\, .$$
Indeed, using~\eqref{equal2c} and equality in~\eqref{cond} for $J$, we have $\displaystyle \sum_{j\in J^c} c_j = \textrm{dim}({\rm span}\{b_j\ ;\ j\in J^c\})$, and also for $\tilde{J} \subset J^c$,  since $P_J^\perp a_j = 0$ for $j\in J$,
\begin{eqnarray*}
\sum_{j\in \tilde{J}}c_j \; + \; \textrm{dim}(V_J) & = & \sum_{j\in J\cup \tilde{J}} c_j \\
& \le & \textrm{dim}(\textrm{span}\{a_j \ ; \ j \in J\cup\tilde{J} \}) \\
& = &  \textrm{dim}(\textrm{span}\{P_J \, a_j\,  + P_J^\perp a_j  \ ; \ j \in J\cup\tilde{J} \} )\\
& \le &  \textrm{dim}(\textrm{span}\{P_J\, a_j \ ; \ j \in J\cup \tilde{J}\}) +  \textrm{dim}(\textrm{span}\{P_J^\perp a_j \ ; \ j \in \tilde{J} \})  \\
 & = & \textrm{dim}(V_J) + \textrm{dim}(\textrm{span}\{b_j \ ; \ j \in \tilde{J} \}) 
\end{eqnarray*}

For an invertible operator $T$ on $\R^n$ we shall use the standard notation 
$$T^{-\ast} := (T^{-1})^\ast= (T^\ast)^{-1}$$
 where $T^\ast x\cdot y = x\cdot Ty$ for all $x,y\in \R^n$. With these definitions, we now state the splitting lemma. Only the first part of the statement is needed to complete the proof of Theorem~\ref{gaussian} ; the rest will be used for the characterization of  extremizers.

\medskip


\begin{lm}\label{peel} Given any family  $A=\{a_1,\ldots, a_m\}$  of $m$ vectors spanning $\R^n$ and $c\in K_A\setminus K_A^\circ$ with $c_j>0$ for all $j\le m$, let $J$ be a 
non-empty proper subset of $\{1,\dots,m\}$  for which equality holds in~(\ref{cond}), and suppose that $J$ has the least cardinality among all such subsets.  Then with $A_J$, $c_J$, $B_{J^c}$ and $c_{J^c}$ defined as above, we have
\begin{equation}\label{decompD}
D(A,c) = D(A_J,c_J) + D(B_{J^c},c_{J^c})\, ,
\end{equation}
and if $D_{\cal G}(B_{J^c},c_{J^c}) = D(B_{J^c},c_{J^c})$, then $D_{\cal G}(A,c) = D(A,c)$.

Suppose next that there exists an extremizing random vector $X$; i.e., a random vector $X$ such that
\begin{equation}\label{equal1}
\sum_{j=1}^m c_j\, S(a_j\cdot X) - S(X) = D(A,c)\ .
\end{equation}
Then
\begin{equation}\label{equal2a}
\R^n = V_J\oplus V_{J^c}\, .
\end{equation}
and this direct sum is an orthogonal decomposition in the inner product given by the covariance matrix of $X$; i.e., $\langle x,y\rangle=   {\mathbb E}\big[(x\cdot (X-\mathbb E X))(y\cdot (X-\mathbb E X))\big] $.

Moreover,  if $T$ is an (invertible) operator on $\R^n$ such that one has the orthogonal decomposition
$$
\R^n = T V_J  \stackrel{\perp}{\oplus}    T V_{J^c}
$$
(for instance $T=H_X^{1/2}$ where $H_X$ is the covariance matrix of  an extremizer $X$, so that
$\langle x,y\rangle = x\cdot H_Xy$), then $X$ is an extremizer~\eqref{equal1} if and only if $T^{-\ast}X$ decomposes as  $T^{-\ast} X =Y+Z$ where $Y$ and $Z$ are \emph{independent} random vectors with values in 
$T V_J$ and $T  V_{J^c}$, and  which are extremizer for $\big([Ta_j \, ; j\in J], c_J\big)$ and 
$\big([Ta_j \, ; j\in J^c], c_{J^c}\big)$, respectively.
\end{lm}


\medskip

The proof of this lemma relies on some well known identities and inequalities concerning conditional entropy  that we now recall. 

Let $E$ and $F$ be two Euclidean spaces (equipped with the Lebesgue measure). If $W$ and $Y$ are two random vectors with values in $E$ and $F$ respectively, with a joint density $\rho(w,y)$ on $E\times F$, let $\rho_Y(y) =  \int_{E} \rho(w,y)\dd w$ and 
$\rho_W(w) =  \int_{F} \rho(w,y)\dd y$ be the two marginal densities on $F$ and $E$, which are of course the densities of $W$ and $Y$ respectively.

Then the {\it conditional density of $W$ given $Y$} is
$\rho(w|y) = \rho (w,y)/\rho_Y(y)$.
The {\em conditional entropy of $W$ given $Y=y$} is then defined to be
$$S(W|Y =y) = \int_{E}\rho(w|y) \ln \rho(w|y) \dd w\ .$$
Since the entropy of $(W,Y)$,  $S(W,Y)$, is given by
$$S(W,Y) = \int_{E\times F} \rho (w,y) \ln \rho (w,y) \dd w\dd y\, ,$$
the identity
\begin{equation}\label{cent1}S(W,Y) = \int_{F}S(W|Y=y) \rho_Y(y)\dd y + S(Y)
\end{equation}
follows directly from the definitions.
Furthermore, by Jensen's inequality
\begin{equation}\label{cent2}
S(W) \le \int_{E}S(W|Y = y)
\rho_Y(y)\dd y\ ,
\end{equation}
and there is equality if and only if $W$ and $Y$ are independent.  

\medskip

\noindent{\bf Proof of Lemma \ref{peel}:} Fix any random vector $X$ with values in $\R^n$ and suppose that $S(X)$ is finite.   
We shall use the definition and notation given before the Lemma. Let $P_J$ denote the orthogonal projection onto $V_J$,  and recall that we have the decomposition~\eqref{equal2c}, so that $P_J^\perp = \id_{\R^n} -P_J$ is also the orthogonal projection onto
${\rm span}\left(\{b_j\ :\ j\in J^c\}\right)$ where $b_{j} = P_J^\perp(a_j)$ for all $j\in J^c$.
Let us introduce
$$Y=P_J X \quad \textrm{ and } Z= P_J^\perp X$$ 
so that $X = Y+Z$. Then $S(X) = S(Y,Z)$ and so from (\ref{cent1}),
\begin{equation}\label{cent3}
S(X) = \int_{V_J}S(Z|Y = y)\rho_Y(y)\dd y + S(Y)\, .
\end{equation}
For each $j\in J$, we have $a_j\cdot X = a_j\cdot Y$, so that
\begin{equation}\label{cent4}
S(a_j\cdot X) = S(a_j\cdot Y) \qquad{\rm for}\quad j \in J\, .
\end{equation}
Note that for $j\in J^c$, $b_{j}\ne 0$, or else $a_j\in V_J$; but this is impossible since $c_j>0$, and 
we already have $\sum_{j=1}^\ell c_j = {\rm dim}(V_J)$.  We have, using the invariance of the entropy under translation,
\begin{equation}\label{ind1}
S(a_j\cdot X | Y = y) =  S(a_j\cdot Z + a_j\cdot y | Y = y)  = S(b_{j}\cdot Z | Y = y) \qquad{\rm for}\quad j \in J^c\, .
\end{equation}
Therefore, by applying~(\ref{cent2}) to $(X\cdot a_j, Y)$ on $\R\times V_J$, we get
\begin{equation}\label{cent5}
S(a_j\cdot X) \le \int_{V_J} S(b_{j}\cdot Z|Y=y)\rho_Y(y)\dd y \qquad{\rm for}\quad j \in J^c\, .
\end{equation}

Now combining (\ref{cent3}),  (\ref{cent4}) and (\ref{cent5}),
we have that
\begin{eqnarray}
\sum_{j=1}^{m} c_jS(a_j\cdot X) - S(X) &\le& \int_{V_J}\left[
\sum_{j\in J^c} c_j S(b_j\cdot Z|Y = y) - S(Z|Y=y) \right] \rho_Y(y)\dd y \nonumber \\
& &  \qquad  + \ \sum_{j\in J} c_j S (a_j\cdot Y) - S(Y)\, \label{cent6}
\end{eqnarray}

It is clear from (\ref{cent6}) and the definition of  $D(B_{J^c},c_{J^c})$ that
$$D(A,c) \le D(A_J,c_J) + D(B_{J^c},c_{J^c}	)\ .$$
To see that there is actually equality here, we use the fact that $J$ 
is a critical set of minimal cardinality.  This implies that  $c_J \in K_{A_J}^\circ$, and by Lemma \ref{interior},
there is a centered Gaussian random vector $Y$ for which
\begin{equation}\label{optimalY}
 \sum_{j\in J}c_j S (a_j\cdot Y) - S(Y) = D(A_J,c_J)\, .
 \end{equation}
Pick $\epsilon>0$ and let  $Z$ be any  random variable with values in $V_J^\perp$ that is independent of $Y$ and
such that
\begin{equation}\label{trial}
\sum_{j\in J^c}c_j S(b_j\cdot Z) - S(Z) > D(B_{J^c},c_{J^c})-\epsilon\ .
\end{equation}
For $\delta > 0$, form  the $\R^n$ valued random vector $X = \delta Y +Z$. 
Since $Y$ and $Z$ are orthogonal and independent, $S(X) = S(\delta Y, Z) = S(\delta Y) + S(Z)$.  The scaling invariance implies that~\eqref{optimalY} holds when $Y$ is replaced by $\delta Y$.
Also, for $j\in J^c$, as $\delta$ approaches zero, $S(a_j\cdot X) = S(b_j\cdot Z + \delta a_j \cdot Y)$ approaches $S(b_j\cdot Z)$.
(Note that by the independence of $Y$ and $Z$, $b_j\cdot Z + \delta a_j Y$ is simply a standard Gaussian regularization of 
$b_j\cdot Z$.)
It now follows that for $\delta$ sufficiently small, 
$$\sum_{j=1}^m c_jS(a_j\cdot X) - S(X) \ge  D(A_J,c_J) + D(B_{J^c},c_{J^c}) -2\epsilon\ .$$
This implies that $D(A,c)\ge D(A_J,c_J) + D(B_{J^c},c_{J^c})$. We have implicitly assumed that $D(B_{J^c},c_{J^c})<+\infty$ (we shall later only need this case, actually), but the argument remains valid if $D(B_{J^c},c_{J^c})=+\infty$. Thus~\eqref{decompD} is established.

Now suppose that $D_{\cal G}(B_{J^c},c_{J^c}) =  D(B_{J^c},c_{J^c})$. Then we may further assume that the random variable $Z$ in the
previous paragraph is a centered Gaussian random variable. Combining this with the independent extremal centered Gaussian random variable $Y$,
provided by Lemma \ref{interior}, we see that we may take the random variable $X$ in the previous paragraph to be
a centered Gaussian. Hence, in this case, $D_{\cal G}(A,c) = D(A,c)$.

It remains to prove the last statements concerning the cases of equality.  

We first assume that we are given a finite entropy random variable $X$ for which  (\ref{equal1}) is satisfied.  By making a translation, we may assume that $X$ is centered; i.e., ${\rm E}(X) = 0$. Furthermore, the covariance matrix is non-degenerate or else the law of $X$ would be concentrated on a proper subspace and this is
inconsistent with finite entropy.  Since $X$ satisfies~(\ref{equal1}), there must be equality in~(\ref{cent6}),
and it must be the case that  
\begin{equation}\label{equal2}
\sum_{j\in J} c_j S(a_j\cdot Y) - S(Y) = D(A_J,c_J)
\end{equation}
and that for each $y\in V_J$,
\begin{equation}\label{equal3}
\sum_{j\in J^c} c_j S(b_j\cdot Z|Y = y) - S(Z|Y=y)  = D(B_{J^c},c_{J^c})\, .
\end{equation}
And since $X$ is centered, so is $Y$.
Next, in addition to equality in~(\ref{equal3}), we must have equality in~(\ref{cent6}). Since the only inequality used in deriving~(\ref{cent6}) was ~(\ref{cent5}),
this in turn requires equality in~(\ref{cent5})
for each $j\in J^c$.  By (\ref{ind1}), this means that for $j\in J^c$,
$$S(a_j\cdot X) = \int_{V_J}S(a_j\cdot X|Y = y)\rho_Y(y)\dd y\ .$$
By the condition for equality in~(\ref{cent2}), this implies that
for $j\in J^c$, $a_j\cdot X$ and $Y$ are independent random variables.  But then for any $y\in V_J$, by independence
$$\langle y,a_j\rangle = {\mathbb E}[(y\cdot Y)(a_j\cdot X)] =  {\mathbb E}(y\cdot Y){\mathbb E}(a_j\cdot X) = 0 \, .$$
This shows that $V_J$ and $V_{J^c}$ are orthogonal subspaces in the inner product defined in terms of the covariance. Thus their dimension sums exactly to $n$ and so~(\ref{equal2a}) holds. 

We now prove the final statement describing how extremizers split. 

Note that given invertible operator $T$ on $\R^n$, a random vector $X$ is extremal~\eqref{equal1} for $(A,c)$ if and only if $T^{-\ast}X$ is extremal for $(TA, c)$ with the notation $TA=[Ta_1, \ldots, Ta_m]$. Indeed, since $a_j\cdot X = Ta_j\cdot T^{-\ast} X$ and $S(T^{-\ast} X) = S(X) + \ln(|\det(T)|)$ we have that~\eqref{equal1} is equivalent to
$$\sum_{j=1}^m c_j\, S(Ta_j\cdot T^{-\ast} X)  - S(T^{-\ast} X) = D(T A, c)$$
and $D(TA, c) = D(A, c) - \ln(|\det(T)|)$.

As in the statement of the lemma, let $T$ be an invertible operator on $\R^n$ such that $\R^n = T V_J  \stackrel{\perp}{\oplus}    T V_{J^c}$. The previous remark explains the mechanism of replacing $A$ by $TA$ and $X$ by $T^{-\ast}X$. So after this transformation we are reduced to proving the statement in the case $T=\id$. Therefore we assume from now on that 
$$\R^n = V_J  \stackrel{\perp}{\oplus}    V_{J^c}.$$
We go back to the beginning of the proof and note  that $b_j = a_j$ for all $j\in J^c$: the orthogonal projection does nothing in this case ($P_{J}^\perp = P_{{J^c}}$).

Assume $X$ is an extremizer~\eqref{equal1} which is decomposed as before as $X=Y+Z$. Then as in the argument above we must have that 
\begin{equation}\label{equal2bis}
\sum_{j\in J} c_j S(a_j\cdot Y) - S(Y) = D(A_J,c_J)
\end{equation}
and that for each $y\in V_J$,
\begin{equation}\label{equal3bis}
\sum_{j\in J^c} c_j S(a_j\cdot Z|Y = y) - S(Z|Y=y)  = D(A_{J^c},c_{J^c})\, ,
\end{equation}
with  $Y$ and $a_j\cdot X$ independent for every $j\in J^c$. Since $a_j\cdot X= a_j\cdot Z$ for every $j\in J^c$ 
we have that $a_j\cdot Z$  is independent of $Y$ for $j\in J^c$ and so $S(a_j\cdot Z|Y = y) =S(a_j\cdot Z)$. Using this together with~\eqref{cent2} for $W=Z$, we get, after integrating~\eqref{equal3bis} with respect to $\rho_Y(y)\, dy$, and applying (\ref{cent2}),
$$
 D(A_{J^c},c_{J^c})\le \sum_{j\in J^c} c_j S(a_j\cdot Z) - S(Z) \, .
$$
By the definition of  $D(A_{J^c},c_{J^c})$ this inequality must be an equality, i.e.
\begin{equation}\label{equal4}
\sum_{j\in J^c} c_j S(a_j\cdot Z) - S(Z)  = D(A_{J^c},c_{J^c})\, ,
\end{equation}
and therefore, there must be equality in the application of ~\eqref{cent2} that we just made.
This implies that $Z$ and $Y$ are independent, as claimed.

Conversely, let $X$ be a random vector such that $X=Y+Z$ in the decomposition $\R^n = V_J  \stackrel{\perp}{\oplus}    V_{J^c}$ with $Y$ and $Z$ independent and such that~\eqref{equal2bis} and~\eqref{equal4} holds. Then we have~\eqref{equal3bis} and we readily check that there is equality at every step. So $X$  is indeed an extremizer~\eqref{equal1}.
\qed


\medskip
We are now ready to prove Theorem \ref{gaussian}:
\medskip

\medskip
\noindent{\bf Proof of Theorem \ref{gaussian}}   By Lemma \ref{interior}, whenever
$c\in K_A^\circ$, $D_{\cal G}(A,c) = D(A,c)$, and there is a Gaussian optimizer. 

Hence it remains to consider the case $c\in K_A\setminus K_A^\circ$. Then taking $J$ to be a proper non-empty subset of $\{1, \ldots, m\}$ of least cardinality for which there is equality in~\eqref{cond}, we may
``peel off'' $|J|$ vectors from our set, as in the first part of Lemma~\ref{peel}, and reduce maters to the consideration of 
$D(B_{J^c},c_{J^c})$. By that Lemma, $D_{\cal G}(A,c) = D(A,c)$ whenever  $D_{\cal G}(B_{J^c},c_{J^c}) = D(B_{J^c},c_{J^c})$.
Now, if $B_{J^c}$ and $c_{J^c}$ are such that for every proper subset of the remaining indices, 
strict inequality holds in the analog of~(\ref{cond}), i.e. $c_{J^c}\in K_{A_{J^c}}^\circ$, then  $D_{\cal G}(B_{J^c},c_{J^c}) = D(B_{J^c},c_{J^c})$ follows from
Lemma~\ref{interior}. Otherwise, we ``peel off'' another proper subset of indices for which equality holds in
(\ref{cond}), and reduce to a problem with a strictly smaller number of vectors. In a finite number of steps, this process must
end.  \qed

\medskip

Our next theorem concerns the cases of equality in the subadditivity inequality. As we have seen in Lemma \ref{peel}, when there is equality,
and no $c_j$ is zero, then
either $c\in K_A^\circ$, or
the variational problem can be split into two problems of the same type, but involving reduced number of vectors, and for random variables taking values in subspaces of a reduced dimension.   

Of course, each of these reduced problems must also have an optimizer, and so
we can apply the same dichotomy to each of them. This leads to the following definition:

\begin{defi}[Totally reducible for $c$)]  \label{def:totred}
Let $A=\{a_1,\dots,a_m\}$ be a family of vectors that spans $\R^n$ and $\{c_1, \ldots, c_m\}$ a set of real numbers with $0\le c_j\le 1$ for $j\le m$. We say that  $\{a_1,\dots,a_m\}$ is \emph{totally reducible for $c$} if $c\in K_A$ and
in case for some $k\ge 1$ there is a decomposition (possibly with $k=1$)
$$\{1,\dots,m\} = J_0 \cup J_1 \cup\ldots \cup J_k$$
where $j\in J_0$ if and only if $c_j=0$, and 
$$\R^n = V_{J_1}\oplus \cdots \oplus   V_{J_k} \quad \textrm {with} \quad V_{J_i} = {\rm span}(\{ a_\ell : \ell \in J_i\})\, ,$$
such  that for each $1\le i \le k$, there is \emph{no} nonempty proper subset of $J_i$
that yields equality in (\ref{cond}). Here, $J_0$ may be empty, but for $1\le i \le k$, $J_i$ is to be non empty. 
\end{defi}

Note that, if  $\{a_1,\dots,a_m\}$ is {totally reducible for $c$}, then we have, with the notation of the definition, that for $1\le k\le m$,
$$c_{J_i}\in K_{J_i}^\circ.$$

The analysis made so far proves the following theorem, which gives a complete analysis of the cases of equality
in the subadditivity inequality. 

\medskip

\begin{thm}\label{gaussian2}
Consider a family   $A = \{a_1,\dots,a_m\}$ of vectors spanning $\R^n$. Then for any
$c\in K_A$, there exists a finite entropy random variable $X$ for which
\begin{equation}\label{equalf}
\sum_{j=1}^m c_jS(a_j\cdot X) - S(X) = D(A,c)\ ,
\end{equation}
if and only if   $\{a_1,\dots,a_m\}$ is totally reducible for $c$. 
In this case, if $\R^n = V_{J_1}\oplus \cdots \oplus   V_{J_k} $ is the corresponding decomposition of $\R^n$ from definition~\ref{def:totred}, let $T$ be any symmetric positive operator on $\R^n$ such that the following orthogonal decomposition holds
\begin{equation}\label{decomS}
\R^n =  T V_{J_1} \stackrel{\perp}{\oplus}  \cdots \stackrel{\perp}{\oplus}   T  V_{J_k} .
\end{equation}
Then the extremizers~\eqref{equalf} are \emph{exactly}  the random vectors $X$ such that $T^{-\ast} X$  decompose as
$$T^{-\ast} X = X_1 + \cdots + X_k$$
where  $\{X_1,\dots,X_k\}$ is an independent set of random variables with each $X_i$ taking values in $T V_{J_i}$ and  extremal for the corresponding problem $(\big([Ta_j \, ; j\in J_i], c_{J_i}\big)$. 
More precisely,
for each $i\le k$, if ${\rm dim}(V_{J_i}) =1$, then $X_i$ can be any finite entropy random variable with values in $TV_{J_i}$;
However, if  ${\rm dim}(V_{J_i}) >1$, then $X_i$ is necessarily Gaussian, and its covariance is a constant multiple of $R_i^2$, where
$R_i$ is the unique positive definite linear transformation on $TV_{J_i}$ such that
$$\sum_{j\in J_i}c_j\left(\frac{R_i Ta_j}{|R_i Ta_j|}\right)\otimes  \left(\frac{R_i Ta_j}{|R_i Ta_j|}\right) = {\rm Id}_{TV_{J_i}} \qquad{\rm and}\qquad
{\rm trace}(R_i^2) = {\rm dim}(T V_{J_i})={\rm dim}(V_{J_i})\ .$$
Finally, if $X$ is an extremizer~\eqref{equalf}  then the symmetric positive operator $T$  defined by $x\cdot T^2 y= {\mathbb E}\big[(x\cdot (X-\mathbb E X))(y\cdot (X-\mathbb E X))\big] $
satisfies the required condition~\eqref{decomS}.
\end{thm}
\medskip

\noindent{\bf Proof:}  The proof relies on successive applications of the Lemmas~\ref{peel} and~\ref{interior}.
First of all, note that the vectors $a_j$ for the indices $j$ such that $c_j=0$ play no role in the inequality, and so  without loss of generality, we may discard these indices without changing $D(A,c)$, the extremizers and $K_A$. So we will assume that $c_j>0$ for all $j\le m$ (this means $J_0=\emptyset$ in the Definition~\ref{def:totred}). 

Assume there exists an extremizer $X$, which, after  translation, can be assumed to be symmetric,  and let $T$ be the symmetric positive operator on $\R^n$ defined by $Tx\cdot Ty = \mathbb E \big[(x\cdot X)(y\cdot X)\big]$. As explained in the proof of the Lemma~\ref{peel} the change of vectors $X\to T^{-\ast}X$ and $a_j\to Ta_j$ reduces the problem to the case $T=\id$, which means that $X$ has unit covariance. Then from Lemma~\ref{peel} we have $\R^n=V_{J_1} \stackrel{\perp}{\oplus} V_{J_1^c}$ for some set of indices $J_1$,
with $c_{J_1}\in K_{A_{J_1}}^\circ$ and $c_{J_1^c}\in K_{A_{J_1^c}}$, and moreover in this orthogonal decomposition $X= X_1+Z$ with $X_1$ and $Z$ extremal for $(A_{J_1},c_{J_1})$ and $(A_{J_1^c}, c_{J_1^c})$, respectively. We apply then Lemma~\ref{peel} on the space $V_{J_1^c}$ to the vector with unit covariance $Z$ which is extremal. This gives for some $J_2\subset J_1^c$ another orthogonal decomposition $V_{J_1^c}= V_{J_2} \stackrel{\perp}{\oplus}V_{J_2^c}$ where  $J_2^c=\{j\in J_1^c\; ; \  j\notin J_2\}$ with $c_{J_2}\in K_{A_{J_2}}^\circ$. After a finite number $k$ of step this process muss end and we have
$$\R^n = V_{J_1}\oplus \cdots \oplus   V_{J_k}$$
with $c_{J_i}\in K_{A_{J_i}}^\circ$ for $i\le k$. This shows that there exists an extremizer only when $\{a_1,\dots,a_m\}$ is totally reducible for $c$. Note that we have also shown that this sum is orthogonal w.r.t. the scalar product given by the covariance of an extremizer.

We assume from now that $\{a_1,\dots,a_m\}$ is totally reducible for $c$ and that $\R^n = V_{J_1}\oplus \cdots \oplus   V_{J_k} $ is the corresponding decomposition of $\R^n$ from definition~\ref{def:totred}.  We can assume that $|J_1|\le |J_2|\le \ldots \le |J_k|$. 
Let $T$ be any symmetric positive operator on $\R^n$ such that the following orthogonal decomposition holds
$$
\R^n =  T V_{J_1} \stackrel{\perp}{\oplus}  \cdots \stackrel{\perp}{\oplus}   T  V_{J_k} .
$$
Of course, there always exists such a linear map $T$. As before the change of vectors $X\to T^{-\ast}X$ and $a_j\to Ta_j$ reduces the problem to the case $T=\id$ and 
$$
\R^n =   V_{J_1} \stackrel{\perp}{\oplus}  \cdots \stackrel{\perp}{\oplus}   V_{J_k} .
$$
With this orthogonal decomposition in hand, we can use Lemma~\ref{peel} to successively ``peel-off" orthogonal blocks. We first apply this Lemma to $J_1$ and $J_1^c=J_2\cup \ldots \cup J_k$, and then on the space $V_{J_1^c}=V_{J_1}^\perp=V_{J_2} \stackrel{\perp}{\oplus}  \cdots \stackrel{\perp}{\oplus}   V_{J_k}$ to $J_2$, and so on.  After $k$ steps we get that $D(A,c)=D(A_{J_1}, c_{J_1}) + \ldots + D(A_{J_k}, c_{J_k})$ and  that a random vector $X$ is an extremizer if and only if it can be written as
$$X=X_1+\ldots +X_k$$
where  $X_i$ has values in $V_{J_i}$ and is extremal for $(A_{J_i}, c_{J_i})$, and with the property that
\begin{equation}\label{indep}
 X_i \quad \textrm{ is independent of }\quad (X_{i+1}, \ldots, X_k) \, , \qquad \textrm{for } i=1, \ldots, k-1 .
\end{equation}
(Note that in order to construct and extremizer $X$ we start with an extremizer $X_k$ on $V_{J_k}$ and, then add an extremal independent $X_{k-1}$ on $V_{J_k}$ in order to get an extremizer on $V_{J_{k-1}}\stackrel{\perp}{\oplus} V_{J_k}$, and so on by repeated applications of Lemma~\ref{peel}).
Observe that the independence  property~\eqref{indep} is equivalent to the independence of the set of random vectors  $\{X_1, \ldots, X_k\}$.
Next remember that for each $i\le m$ we have $c_{J_i} \in K_{J_i}^\circ$. Thus Lemma~\ref{interior} applies and when ${\rm dim}(V_{J_i}) >1$ then $X_i$ is Gaussian and its variance is imposed as stated. Recall that in dimension $1$ the problem is trivial and all random variables are extremal (in particular Gaussian variables are extremal).

\qed

\medskip

Note that the previous theorem tells in particular that when there exists optimizers, there exists Gaussian optimizers (however this was not a needed step in our approach).

Of course, by Theorems~\ref{dual} and~\ref{cases}, we now also know that optimizers for the classical Brascamp--Lieb inequality
exist under the exact same conditions for optimality described in Theorem~\ref{gaussian2}, and that moreover, the optimizers
Brascamp--Lieb inequality are exactly the marginals of the optimizing probability densities for the subadditivity inequality. The full description of optimizers (in one dimensional Brascamp-Lieb inequalities) was given in~\cite{CLL1}, building on a previous characterization by Barthe~\cite{B2}. In the multidimensional case, building on Barthe's work too, Bennett-Carbery-Christ-Tao~\cite{BCCT1}  obtained some description, but the problem was completely solved only recently by Valdimarsson~\cite{V}.

\medskip


\section{Consequences of the general subadditivity inequality in $\R^n$}

There are several interesting consequences of Theorems \ref{gaussian} and \ref{gaussian2}.  The first  is a generalization of Hadamard's inequality for determinants:

\medskip
\begin{thm}\label{had} Consider any family  $A=\{a_1,\dots,a_m\}$ of $m$  vectors that span $\R^n$, any set of numbers
$\{c_1,\dots,c_m\}$ with $0\le c_i \le 1$. Then with $D(A,c)$ as above, for any linear transformation $T$ from $\R^n$ to $\R^n$,
\begin{equation}\label{had1}
|\det(T)| \le e^{D(A,c)}\left(\prod_{j=1}^m |T(a_j)|^{c_j}\right)\ ,
\end{equation}
and this inequality is sharp in that the constant $e^{D(A,c)}$ cannot be decreased. Moreover, for $c\in K_A^\circ$, there is transformation $T$
with $\det(T) =1$ for which equality holds in (\ref{had1}),and, when $n\ge 2$, if we take $T$ to be positive, then $T$ is unique (up to multiplication by a positive scalar).
\end{thm}
\medskip

\noindent{\bf Remark:} In the case that $m=n$, and the vectors $\{a_1,\dots, a_m\}$ are an orthonormal basis, and $c_1 = \cdots = c_n = 1$,
this reduces Hadamard's inequality for determinants. In the special case $\sum_{m=1}^mc_ja_j\otimes a_j = \id_{\R^n}$, so that $D(A,c)= 0$, 
this result has been proved by Ball \cite{Ball}, with a very simple proof.

For simplicity we have stated the existence of an extremal $T$ only when $c\in K_A^\circ$, but the right condition is that $A$ is totally reducible for $c$, just as in Theorem~\ref{gaussian2}.

\medskip

\noindent{\bf Proof:}  By making a polar decomposition, we may assume without loss of generality that $T$
 is positive definite.  Let $G_T$ be the centered Gaussian random variable with ${\rm E}(u\cdot G_T) = |T(u)|^2$ for all vectors $u$ in
 $\R^n$.  Then simply evaluating  the left hand side of $\sum_{j=1}^m c_jS(a\cdot G_T) - S(G_T) \le D(A,c)$, we obtain 
 (\ref{had1}).  Then Theorem \ref{gaussian} provides the rest.  \qed
 
 \medskip

Theorem \ref{had} gives us one simple variational expression for $D(A,c)$, namely
$$D(A,c) = \sup\left\{ \ln\left(\frac{ |\det(T)|}{\prod_{j=1}^m |T(a_j)|^{c_j}}\right)\ :\ T\ {\rm positive\  definite}\ \right\}\ .$$
There is however a simpler variational formula for $D(A,c)$ over an even lower dimensional space, as suggested by the fact that  $e^{D(A,c)}$ is also the sharp constant in the Brascamp--Lieb inequality. By the classical theorem of Brascamp and Lieb, 
$e^{D(A,c)}$ may be computing by taking the functions $\{f_1,\dots,f_m\}$ in the  Brascamp--Lieb inequality to be centered Gaussians; i.e.,
$$\{f_1(t),\dots,f_m(t)\} = \{e^{-(s_1 t)^2}, \dots, e^{-(s_m t)^2}\}\ ,$$
and varying the $m$ numbers $s_1,\dots,s_m$. 
This leads directly to the variational expression (\ref{gp1}) for $D(A,c)$. Let us recall that the existence of optimizers for this problems
was proved by Brascamp and Lieb \cite{BL} under the hypothesis that every set of $n$ vectors chosen from $\{a_1,\dots,a_m\}$ is linearly
independent and later proved by Barthe \cite{B2} for $c\in K_A^\circ$.  The next theorem gives the complete result.  Although
the variational formula  (\ref{gp1}) can be deduced by duality, we give a direct proof of it starting from the subadditivity inequality.

\medskip
\begin{thm}\label{gp} Consider any set $\{a_1,\dots,a_m\}$ of $m$  vectors that span $\R^n$, $n\ge 2$. Let 
$\{c_1,\dots,c_m\}$ be any set of numbers with $0\le c_i \le 1$ verifying~\eqref{sumcon}. 
Let $T$ denote the $m\times m$ diagonal matrix whose $j$th diagonal entry is $t_j$, and  define the function $\Phi_A(t_1,\dots,t_m)$ by
$$\Phi_A(t_1,\dots,t_m) =   \ln\det (A e^T A^t)\ .$$
This is a convex function on $\R^m$, and
\begin{equation}\label{gp1}
D(A,c) + \sum_{j=1}^m c_j\ln(c_j) = 
\frac{1}{2}\sup_{\{t_1,\dots,t_m\}}\left( \sum_{j=1}^m c_jt_j -  \Phi_A(t_1,\dots,t_m) \right) \ .
\end{equation}
The supremum in (\ref{gp1}) is attained if and only if  $\{a_1,\dots,a_m\}$ is totally reducible for $c$.
Moreover, 
\begin{equation}\label{gp2}
\Phi_A(t_1,\dots,t_m) =  \sup_{\{c_1,\dots,c_m\}}\left( \sum_{j=1}^m c_jt_j - 2\left(D(A,c) + \sum_{j=1}^m c_j\ln(c_j)\right)\right)\ .
\end{equation}
\end{thm}
\medskip

\noindent{\bf Proof:}  For an $m\times m$ diagonal  matrix $S$ with positive entries $s_j$, introduce
 $R_S := ((AS)(AS)^t)^{-1/2}$ as in the proof of  Lemma~\ref{fjcond}. Let $G$ be a standard Gaussian random vector on $\R^n$ (i.e., $G \in \mathcal N (0, \id)$),  and set $G_S= R_S G$. Then
$$\sum_{j=1}^mc_j S(a_j\cdot G_S)  - S(G_S) =  -\ln(\det(R_S^{-1}))  - \frac{1}{2}\sum_{j=1}^mc_j\ln(|R_Sa_j|^2)\ .$$
However,
$$|R_Sa_j|^2 = s_j^{-2}|R_S (s_ja_j)|^2=   s_j^{-2}|R_S(SA)e_j|^2 = s_j^{-2}e_j\cdot(AS)^t((AS)(AS)^t)^{-1}(AS)e_j\ ,$$
where $e_j$ is the $j$th standard basis vector in $\R^m$.  Recall that  
$e_j\cdot(AS)^t((AS)(AS)^t)^{-1}(AS)e_j$   is the $j$th diagonal entry of the orthogonal projection in $\R^m$ onto the image of $(AS)^t$.  Since this  orthogonal projection has rank $n$,
its trace is $n$. Therefore, if we define $c_j(S) =  |R_S(SA)e_j|^2$, we have
$$\sum_{j=1}^m c_j(S) = n$$
for all $S$.  Thus, by Jensen's inequality,
$$\sum_{j=1}^m c_j\ln(c_j) \ge  \sum_{j=1}^m c_j\ln(c_j(S))\ ,$$
with equality exactly when $c_j(S) = c_j$ for all $j$.  Therefore, for all $S$,
$$D(A,c) \ge \sum_{j=1}^mc_j S(a_j\cdot G_S)  - S(G_S) \ge  - \ln(\det(R_S^{-1}))    + \frac{1}{2}\sum_{j=1}^mc_j\ln(s_j^2)  - \frac{1}{2}\sum_{j=1}^m c_j\ln(c_j)\ .$$
so that 
$$D(A,c) + \sum_{j=1}^m c_j\ln(c_j) \ge  \frac{1}{2}\left( \sum_{j=1}^mc_j\ln(s_j^2)  - \ln \det(R_S^{-2} ) \right)\ .$$
Moreover, as we see from  the proof of Lemma~\ref{fjcond} (based on an observation by Barthe) and Lemma~\ref{interior} and the remarks made just above, there is equality when $c\in K_A^\circ$ and  $S = S_0$ is  the choice of $S$
(unique up to a multiple) 
for which (\ref{el}) is true. Let $T$ denote the $m\times m$ diagonal matrix whose $j$th diagonal entry is $t_j = \ln s_j^2$. Then
$\ln(\det(R_{S}^{-2})) = \ln(\det (A e^T A^t)$ and therefore, if we define the function $\Phi_A$ by
$$\Phi_A(t_1,\dots,t_m) =  \ln \det (A e^T A^t)\ ,$$
we have, for every $t_1, \ldots, t_m \in \R$
\begin{equation}\label{dualityPhi}
2 \Big(D(A,c) + \sum_{j=1}^m c_j\ln(c_j) \Big) + \Phi_A (t_1, \ldots, t_m) \ge \sum_{j=1}^m c_j \, t_j
\end{equation}
with equality, when $c\in K_A^\circ$ for some choice of $t_j$'s. The function $c\longrightarrow 2 D(A,c) + 2 \sum_{j=1}^m c_j\ln(c_j) $ is convex (because, as mentioned at the beginning of the previous section, the function $c\to D(a,c)$ is convex by definition), and its domain (i.e. where it is $<+\infty$) is $K_A$. Therefore we get that
$$
\Phi_A (t_1,\dots,t_m) = \sup_{c\in K_A^\circ } \Big\{  \sum_{j=1}^m c_j \, t_j - 2 \Big(D(A,c) + \sum_{j=1}^m c_j\ln(c_j) \Big) \Big\} =   \sup_{c\in K_A} \Big\{  \ldots \Big\} =  \sup_{c\in \R^m} \Big\{  \ldots \Big\}.
$$
This shows that $\Phi_A$ is convex on $\R^m$ and that it is the Legendre transform of the convex function  $c\longrightarrow 2 D(A,c) + 2 \sum_{j=1}^m c_j\ln(c_j) $.

Moreover,  for given $A$ and $c$, equality in~\eqref{dualityPhi}  for some $t_1,\dots,t_m$
means that for the corresponding values $s_1,\dots,s_m$, the Gaussian $G_S$ is an extremizer for
the variational problem defining $D(A,c)$. By Theorem \ref{gaussian2}, tis means that 
$\{a_1,\dots,a_m\}$ is totally reducible for $c$.

Conversely, if $\{a_1,\dots,a_m\}$ is totally reducible for $c$, then the variational problem 
in  (\ref{gp1}) splits into a sum of independent and orthogonal (after a suitable linear transformation $T$) such problems, but of the interior type 
(i.e. $c\in K_{TA}^\circ$) for which Barthe showed optimiziers to exist. 
Equivalently, the next Theorem~\ref{fjcond2}  ensures that we can find a positive operator $R$ for which the decomposition of the identity~\eqref{propR} holds. Then, as mentioned in the remark after the proof of Lemma~\ref{interior}, the random vector $RG$ is extremal for $D(A,c)$ and setting $s_j^2=c_j/ |Ra_j|^2$ we have that $R=R_S$ and $c_j(S)=c_j$ by construction (see the proof of Lemma~\ref{fjcond}). This guaranties equality at all steps of our computation above and thus ensures equality in~\eqref{dualityPhi}
\qed

\medskip

\noindent{\bf Remark:}  We have proved that
 $$D(A,c) + \sum_{j=1}^m c_j\ln(c_j) = \frac{1}{2}\Phi_A^*(c)\ $$
where $\Phi_A^*$ denotes the Legendre transform of $\Phi_A$.
Since $\nabla \Phi_A^*(\nabla \Phi_A(0)) = 0$, the choice $c= \nabla \Phi_A(0)$ minimizes 
$\Phi_A^*(c)$, and hence $D(A,c) + \sum_{j=1}^m c_j\ln(c_j)$.  There is a misprint in \cite{CLL1} in which it is stated (in slightly different notation) that this choice of
$c$ minimizes $D(A,c)$ itself.

\medskip

We finally return to Lemma ~\ref{fjcond}, as we are now in a position to give necessary and sufficient conditions for the existence
of the change of variables provided there. 

Let $A = \{a_1,\dots, a_m\}$ be family of $m$ vectors spanning $\R^n$, and let $c$ be any vector in $\R^m$ with
$0\le c_j \le 1$ for all $j$.   Theorem \ref{gaussian2} gives us necessary and sufficient conditions for the existence
of an extremal $X$ for the subadditivity inequality. By Theorem \ref{cases}, these conditions  are also 
necessary and sufficient for the existence of extremals for the Brascamp-Lieb inequality.  Moreover, we see that extremals for the latter exist if and only if centered Gaussian extremals exist. 

From here, it is easy to prove the following theorem which supersedes Lemma \ref{fjcond}, and gives necessary and sufficient
conditions for the existence of the change of variables considered there.  This result was obtained (in the more general multidimensional setting)  by  Bennett-Carbery-Christ-Tao~\cite{BCCT1} along their study of the Brascamp-Lieb extremizers ; here we use the extremizers to the subadditivity of entropy inequality. Though this theorem concerns a problem in linear algebra, we do not know a direct proof of it in a purely linear algebra context, though there may be one.

\medskip
\begin{thm}\label{fjcond2} Let $\{a_1,\dots,a_m\}$ be any collection of $m$ vectors that span $\R^n$ for $n\ge 2$.
Let $\{c_1,\dots,c_m\}$ be any $m$ numbers satisfying $0 \le c_j \le 1$ for each $j$.  Then  there exists an
an invertible symmetric matrix $n\times n$ matrix $R$ so that
$$\sum_{i=1}^m c_j\left(\frac{Ra_j}{|Ra_j|}\right)\otimes \left(\frac{Ra_j}{|Ra_j|}\right)= \id_{\R^n}\, $$
if and only if  the set $\{a_1,\dots,a_m\}$ is totally reducible for $c$
\end{thm}

\medskip

\noindent{\bf Proof:}  The proof of Lemma \ref{interior} shows that whenever such a matrix $R$ exists, there exists an optimizer for
the subadditivity inequality. Thus, by Theorem \ref{gaussian2}, the condition that $\{a_1,\dots,a_m\}$ is totally reducible for $c$
is necessary. 

Conversely assume that  $\{a_1,\dots,a_m\}$ is totally reducible for $c$ 
 and that $\R^n = V_{J_1}\oplus \cdots \oplus   V_{J_k} $ is the corresponding decomposition of $\R^n$ from definition~\ref{def:totred}.   We can then find an invertible operator  $T$ on $\R^n$ such that the following orthogonal decomposition holds
$$
\R^n =  T V_{J_1} \stackrel{\perp}{\oplus}  \cdots \stackrel{\perp}{\oplus}   T  V_{J_k} .
$$
Since we have $c_{J_i}\in K_{A_{J_i}}^\circ = K_{TA _{J_i}}^\circ$ for $i\le m$ (with $TA_{J_i}=[Ta_j, j\in J_i]$), we may use the Lemma \ref{fjcond} on each of the reduced orthogonal subspaces $TV_{J_i}$ ; this gives us some symmetric  invertible operator $R_i$ on $TV_{J_i}$ and putting all  these operators together we get a symmetric invertible  operator $R$ on $\R^n$ such that
$$\sum_{i=1}^m c_j\left(\frac{RTa_j}{|RTa_j|}\right)\otimes \left(\frac{RTa_j}{|RTa_j|}\right)= \id_{\R^n}\, 
$$
Then the positive symmetric operator $\tilde{R}=\sqrt{T^{\ast} R^2T}$ satisfies the desired property
$$\sum_{i=1}^m c_j\left(\frac{\tilde{R}a_j}{|\tilde{R}a_j|}\right)\otimes \left(\frac{\tilde{R}a_j}{|\tilde{R}a_j|}\right)= \id_{\R^n}\, .$$

\qed

\medskip


\section{A convolution inequality for eigenvalues}

We investigate here the dual of the superadditivity of Fisher information inequality~\eqref{infosub} from Proposition~\ref{fisher}.

In Section  2 we have shown that the Legendre transform of the entropy provides an equivalence between subadditivity of the entropy and Brascamp-Lieb inequalities. It turns out that the Fisher information is also a convex functional and its Legendre transform is known to be
the smallest eigenvalue of a Schr\"odinger operator. (This is used extensively  in the theory of large deviations, for example). We shall use this fact to derive  a subadditivity of the smallest eigenvalues of Schr\"odinger operators.

\medskip

For any continuous bounded function $V$ on $\R^n$, define  
\begin{equation}\label{rayritz}
\lambda(V) = \sup\left\{  \int_{\R^n} V(x)  \phi^2(x)\dd x -  4\int_{\R^n}|\nabla \phi(x)|^2\ : \ \int_{\R^n}\phi^2(x)\dd x =1\ \right\}\ 
\end{equation}
Then $-\lambda(V)$ is the ``ground state''  eigenvalue of
$$-4\Delta - V\ ,$$
provided the bottom of the spectrum is an eigenvalue, and in any case, it is the bottom of the spectrum.

Then since 
$$I(f) = \int_{\R^n}\frac{|\nabla f|^2}{f}\dd x = 4\int_{\R^n}|\nabla \sqrt{f}|^2\dd x\ ,$$
we can rewrite (\ref{rayritz}) as
$$\lambda(V) = \sup\left\{  \int_{\R^n} V(x)  f(x)\dd x -  I(f)\ \right\}\ ,$$
where the supremum is taken over all probability densities $f$. This gives us the analog of (\ref{ful}) for Fisher information:
\begin{equation}\label{ful2} \int_{\R^n} V(x)  f(x)\dd x  \le \lambda(V) + I(f)\ ,
\end{equation}
with equality if and only if $f = \phi^2$ where $(-4\Delta - V)\phi = -\lambda(V)\phi$. 
(Here, by the definition (\ref{rayritz}) of $\lambda(V)$, $\phi$ is the ``ground state'' eigenfunction.

Now let $V_1, \dots V_n$ be continuous functions on $\R$, and define
$$V = \sum_{j=1}^n V_j(u_j\cdot x)\ $$
where $\{u_1,\dots, u_n\}$ is any orthonormal basis for $\R^n$.  Then 
$$-4\Delta - V = \sum_{j=1}^n\left(-4(u_j\cdot \nabla)^2 - V_j(u_j\cdot x)\right)\ ,$$
so that, by separation of variables,
$$\lambda(V) = \sum_{j=1}^n\lambda(V_j)\ .$$
The following result generalizes this to the case in which we have $m$ unit vectors $\{u_1,\dots,u_m\}$ 
satisfying (\ref{fj}):

\medskip

\begin{thm}\label{ev}
Let $\{u_1,\dots, u_m\}$ be any $m$ \emph{unit} vectors in $\R^n$ such that there are positive numbers
$c_1,\dots, c_m$ satisfying
$$\sum_{j=1}^m c_j\, u_j\otimes u_j = \id_{\R^n}\, .$$
For any $m$ continuous bounded functions $V_1,\dots,V_m$ on $\R$, define on $\R^n$
$$V(x) = \sum_{j=1}^m V(u_j\cdot x)\ .$$
Then
\begin{equation}\label{ev2}
\lambda(V) \le   \sum_{j=1}^m c_j\lambda\left( \frac{1}{c_j}V_j\right)\ .
\end{equation}
\end{thm}

\medskip
\noindent{\bf Proof:}  Choose an $\epsilon>0$ and a probability
density $f = \phi^2$ 
such that
$$ \int_{\R^n} V(x)  \phi^2(x)\dd x -  4\int_{\R^n}|\nabla \phi(x)|^2 \ge \lambda(V) -\epsilon\ .$$
Then  using~\eqref{infosub},
\begin{eqnarray*}
\lambda\left(V \right)-\epsilon &\le&  \int_{\R^n} f(x)\left(\sum_{j=1}^m V_j(u_j\cdot x)\right)\dd x - I(f)\nonumber\\
&=& \sum_{j=1}^m \int_{\R} f_{(u_j)}(t) V_j(t)\dd t - I(f)\nonumber\\
&\le& \sum_{j=1}^m \int_{\R} f_{(u_j)}(t) V_j(t)\dd t  -  \sum_{j=1}^m c_jI(f_{(u_j)})\nonumber\\
&=&  \sum_{j=1}^m c_j\left(\int_{\R} f_{(u_j)}(t) \left(\frac{1}{c_j}V_j\right)(t)\dd t  -  I(f_{(u_j)}) \right)\nonumber\\
&\le&\sum_{j=1}^m c_j \lambda \left(\frac{1}{c_j}V_j\right)\ , 
\end{eqnarray*} 
Since $\epsilon>0$ is arbitrary, this proves the result.
\qed

\medskip

The inequality (\ref{ev2}) is sharp since one can use another Legendre transform, as in the proof of Theorem \ref{dual}, and see that it implies the sharp inequality (\ref{infosub}).
Inequality~\eqref{ev2} could also be proved using a semi-group (or Stochastic) method inspired by  the one used by Borell~\cite{borell00} in his study of  Brunn-Minkowski type inequalities (which, somehow,  are the \emph{converse} of the inequalities considered here);   this  would be more complicated  than starting from the inequality~\eqref{infosub} for the Fisher information, though. 

An analogous result for functions on the sphere could be given using the sharp superadditivity
of Fisher information inequality proved in~\cite{BCM}.

\end{document}